\documentclass[twocolumn]{autart}    % Enable this line and disable the 
\usepackage{cite}
\usepackage{har2nat}
\usepackage{natbib}
\usepackage{url}
\usepackage{amsmath}
\allowdisplaybreaks
\usepackage{amssymb}
\usepackage{enumitem}
\usepackage{bm}
\usepackage{algorithm,algorithmic}
\usepackage{mathtools}
\usepackage{makecell}
\mathtoolsset{showonlyrefs}
\usepackage{refcount}
\usepackage{subfigure}
\usepackage{wrapfig}
\usepackage{footnote}
\usepackage{threeparttable}
\usepackage{color}\definecolor{RED}{rgb}{1,0,0}\definecolor{BLUE}{rgb}{0,0,1}
\newtheorem{exam}{Example}

\DeclareMathOperator*{\minimize}{minimize}

\DeclareMathOperator*\argmax{\text{argmax}}

\def\cI{\mathcal{I}}

\def\ln{\text{ln}}

\def\tG{\overline{G}}
\def\brd{\textrm{Broyd}}
\def\gbd{\textrm{gBroyd}}
\def\R{\mathbb{R}}
\def\Rn{\mathbb{R}^{n}}
\def\Rnn{\mathbb{R}^{n\times n}}
\def\tp{\mathsf{T}}
\def\df{\textrm{def}}
\def\equdef{\overset{\df}{=}}

\def\bea{\begin{equation}\begin{alignedat}{-1}}
\def\ena{\end{alignedat}\end{equation}}
\def\bee{\begin{equation}}
\def\ene{\end{equation}}

\def\bigO{\mathcal{O}}

\def\T{\mathsf{T}}

\begin{document}

\begin{frontmatter}

\title{Distributed Adaptive Greedy Quasi-Newton Methods with Explicit Non-asymptotic Convergence Bounds} % Title, preferably not more                                % than 10 words.

% \thanks[footnoteinfo]{This paper was not presented at any IFAC meeting.}
\author{Yubo~Du},
\ead{duyb19@mails.tsinghua.edu.cn}
\author{Keyou~You\corauthref{cor}}
\ead{youky@tsinghua.edu.cn}

\corauth[cor]{Corresponding author}
\address{Department of Automation and BNRist, Tsinghua University, Beijing 100084, China}                                                    % full addresses

                             % keyword list or with the 
                                          % help of the Automatica 
                                          % keyword wizard

\begin{abstract}                          % Abstract of not more than 200 words.
Though quasi-Newton methods have been extensively studied in the literature, they either suffer from local convergence or use a series of line searches for global convergence which is not acceptable in the distributed setting. In this work, we first propose a {\em line search free} greedy quasi-Newton (GQN) method with adaptive steps and establish explicit non-asymptotic bounds for both the global convergence rate and local superlinear rate. Our novel idea lies in the design of multiple greedy quasi-Newton updates, which involves computing Hessian-vector products,  to control the Hessian approximation error, and a simple mechanism to adjust stepsizes to ensure the objective function improvement per iterate. Then, we extend it to the master-worker framework and propose a distributed  adaptive GQN method whose communication cost is comparable with that of first-order methods, yet it retains the superb convergence property of its centralized counterpart. Finally, we demonstrate the advantages of our methods via numerical experiments.  
\end{abstract}

\end{frontmatter}

\section{Introduction}\label{sec:intro}
Quasi-Newton methods contain a popular class of optimization algorithms where the exact Hessian in the Newton method is iteratively approximated, and  have been extensively studied for quite a long history  \citep{nocedal2006numerical}.  Though they retain the superlinear convergence property of Newton  methods near the optimal solution, they either suffer from local convergence with a small domain or require a series of line searches for global convergence whose number of searches per iterate is unpredictable.  

Recently,  \citet{polyak2020new} proposes a novel procedure for directly adjusting stepsizes in Newton methods to achieve a wide range of convergence (even globally) and local quadratic convergence.  Unfortunately, such an adaptive idea cannot be directly applied to the quasi-Newton methods (e.g., DFP, BFGS, SR1) as their Hessian approximation errors cannot be evaluated or controlled. In this work, we adopt the greedy quasi-Newton (GQN) update in \citet{rodomanov2021greedy} which is a variant of  the Broyden family with the greedily selected vectors to maximize a certain measure of progress and involves computing Hessian-vector products (HVPs). To ensure global convergence, we propose to use {\em multiple} GQN updates to uniformly bound the Hessian approximation errors, and then design a new simple mechanism to adaptively adjust stepsizes, leading to the line search free AGQN of this work. Under the strongly convex objective function with Lipschitz continuous Hessian, we derive explicit non-asymptotic  bounds  for both the global convergence rate and local superlinear rate. It is worthy mentioning that even for quasi-Newton methods, the (local) non-asymptotic convergence has been unknown until very recently \citep{rodomanov2022rates,jin2022non,rodomanov2021new}. Our AGQN is the {\em first} QN method that is line search free and achieves global convergence with local superlinear rates. 

Next,  we extend to the distributed master-worker setting where the master updates the decision vector by using information from all workers. Though the distributed method can significantly accelerate the centralized one, how to fully explore the  potential of second-order methods in the distributed setting largely remains open. One of the major bottlenecks lies in that 
 any worker cannot  compute the QN direction by itself and has to use information from other workers.  This is in sharp contrast to  the first-order methods where each worker independently computes its gradient direction in parallel, and thus it is relatively easy to achieve speed-up property.  
 
 A natural idea is that the worker only transmits a limited number of vectors  to form the Hessian estimate or the QN direction. We propose a divide-and-combine strategy to ensure that each worker only sends a vector of dimension $O((\tau+1)n)$ per iterate to the master where $n$ is the dimension of decision vector and the integer valued $\tau\ll n$ is newly introduced to ensure global convergence. 
 Noticeably, $\tau$ can be as small as one, and even zero if the algorithm enters the local region of superlinear convergence, which can be explicit in practice. From this perspective, the communication complexity per iterate is comparable to that of the first-order method in the master-worker setting.  Yet our distributed method, namely DAGQN, not only  retains the superb convergence property of the AGQN with non-asymptotic convergence bounds, but also is line search free.  Specifically, we first divide the task of approximating Hessian among all workers, i.e., let each worker $i$ perform GQN updates to approximate its individual local Hessian and extracts the required vector for the GQN update  to transmit to the master in a single packet. Since the master can reconstruct the same copies of local Hessian estimates, it simply aggregates them to compute the QN direction. Besides the low communication complexity, we are still able to establish explicit non-asymptotic bounds for both the global convergence rate and local superlinear rate for our distributed GQN method. 
 
Many efforts have been made to reduce the communication cost for the distributed setting.  For example,  each worker in NL \citep{Islamov2021} and FedNL \citep{Safaryan2022} dynamically learns the Hessian at the optimum  and sends compressed shift term to the master. \citet{Agafonov2022} builds upon FedNL with a more computation and memory efficient scheme. GIANT \citep{Wang2018}  learns the Newton direction via the harmonic mean of local Newton directions from workers, and each worker in DANE \citep{Shamir2014} sends a solution vector of a local subproblem  to the master. DINGO \citep{Crane2019} and DINO \citep{crane2020dino} find descent directions by combining local Newton directions from workers. \citet{Fabbro2022} approximates Hessian with its principal eigenvectors, and transmits them in an incremental and sporadic manner to amortize the communication cost. \citet{Alimisis2021} vectorizes Hessian and performs lattice-based quantization before communication. DAve-QN \citep{Soori2020} approximates Hessian with the classical BFGS but implements it using an asynchronous scheme. In comparison, our DAGQN has at least two sharp differences: (a) it only involves computing HVPs for updating, instead of explicitly computing Hessian or solving a subproblem per iteration.  (b) The global convergence is achieved via adaptively tuning stepsize without line search  and the local convergence has explicit superlinear rates.  
  
Note that there are also distributed second-order methods for peer-to-peer networks, e.g., DAN-LA \citep{zhang2022distributed} estimates the Hessian via an outer product of a moderate number of vectors that are obtained via the low-rank matrix approximation. And the decentralized BFGS method in \citet{eisen2017decentralized} has a linear rate. See \citet{Mokhtari2016,eisen2017decentralized,lee2018distributed,Tutunov2019,Soori2020, Mansoori2020, Li2021, Li2022,zhang2022distributed}.

 Finally, we adopt a logistic regression problem with real datasets to validate the advantages of our proposed methods.  Numerical results validate that both AGQN and its distributed version are quite competitive for problems where HVPs can be efficiently computed. 
Overall, contributions of this work can be summarized as follows:
\begin{enumerate}[label=(\alph*)]
	\item Our AGQN is the first QN method that is not only line search free, but also achieves global convergence and local superlinear convergence with explicit non-asymptotic bounds.
	\item  The AGQN is extended to the master-worker setting with a low communication complexity that is comparable to that of distributed first-order methods.   Yet, it retains the superb convergence property of the AGQN. 
	\end{enumerate}

The rest of this paper is organized as follows. In Section \ref{sec_gqn}, we briefly review the GQN method. In Section \ref{sec1_global_convergence}, we propose the AGQN method and study its convergence. It should be noted that the preliminary results of this section have been presented in \citet{du2022adaptive}.   In Section \ref{sec:GQN-LA} and \ref{sec:converge}, we propose the distributed AGQN method and conduct convergence analysis, respectively. In Section \ref{sec:experiments}, we perform simulations to validate the performance of our methods.

{\bf Notation:}  $\Vert \cdot \Vert$ is the standard Euclidean norm, and $\langle\cdot,\cdot \rangle$ is the standard inner product in $\Rn$. For a twice differentiable function $f: \Rn \rightarrow  \R$, denote its gradient and Hessian as  $\nabla f$ and $\nabla^{2} f$, respectively.  For a symmetric matrix $A$, denote $A \succ 0$ $(A  \succeq 0)$ if $A$ is positive (semi-)definite. Given $A \succeq 0$ and $B \succeq 0$, the relation $A \succeq B$ means that $A-B \succeq 0$ and $
\langle A, B \rangle = \textrm{Trace}(AB).$ 
For simplicity, we write 
$\Vert  h\Vert_{A} = \langle Ah, h\rangle^{1/2}
$ 
and 
$\Vert h \Vert_{A}^{*} = \langle A^{-1}h, h \rangle^{1/2}, \forall h\in \Rn.
$
For simplicity, we omit the iteration index $k$ and use $x$, $x^-$, $x^+$ to represent  $x^k$, $x^{k-1}$, $x^{k+1}$, respectively. Similar notation will be applied to other vectors and variables.

\section{The basics of greedy quasi-Newton and our objective }
\label{sec_gqn}
In this section, we review some basics of the greedy quasi-Newton (GQN) method in \citet{rodomanov2021greedy}.  In the literature, three QN methods have been widely used, see e.g., \citet{nocedal2006numerical}.  Given symmetric positive definite matrices  $A,G \in \Rnn$ and a vector $u\in\R^n$, satisfying that $G \succeq A$ and $Gu \neq Au$, the  SR1 update is
\bea\label{equ:SR1_update}
\textrm{SR1}(G, A, u) = G - \frac{ (G - A) uu^{\T} (G - A)}{\langle (G - A)u, u \rangle}
\ena
where $u\in\Rn$, and the DFP update is given by
\bea\label{equ:DFP_update}
\textrm{DFP}(G, A, u) =~& G - \frac{Auu^{\T}G + Guu^{\T}A}{\langle Au, u \rangle} \\
& + \left( \frac{\langle Gu, u \rangle}{\langle Au, u \rangle} + 1\right) \frac{Auu^{\T}A}{\langle Au, u \rangle}.
\ena

A convex combination of the SR1  and the DFP leads to the celebrated \textit{Broyden family}  
\bea\label{equ:broyden_update_formula}
&\brd(G, A, u) =  \\
&~~~~~~~~~~~~\phi\textrm{DFP}(G, A, u)
+ (1 - \phi)\textrm{SR1}(G, A, u)
\ena
where $\phi$ can be chosen any value from $[0,1]$ at every call of \eqref{equ:broyden_update_formula}.  Let $\phi  = {\langle Au, u \rangle}{\langle Gu, u \rangle}^{-1} \in (0, 1]$, it reduces to the  BFGS  update
\bea\label{equ:BFGS_update}
\textrm{BFGS}(G, A, u) = G - \frac{Guu^{\T}G}{\langle Gu, u \rangle} + \frac{Auu^{\T}A}{\langle Au, u \rangle}.
\ena
In the optimization progress, we use \textsf{Broyd} to update the estimate of the Hessian matrix  via two vectors $u$ and $Au$, and thus avoid directly computing the Hessian matrix. If $A = \nabla^{2}f$, then the {Hessian-vector} product (HVP) $Au$ can be computed as a whole with  similar computational efforts as that of $\nabla f$ \cite{pearlmutter1994fast}. We use the following example to illustrate this point. 

\begin{exam}\label{example:1}
	In applications, we usually require to solve empirical risk problems (ERP) in the form of  the following generalized linear model (GLM)
	\bea
	\varphi(x) = \sum_{i=1}^{m}\varphi_{i}(a_{i}^{\T}x)
	\ena
	where $x\in \Rn$, $\varphi_{i}$ is a twice differentiable function and $a_{i}$ is a training sample. In the GLM, it is trivial that
	\bea
	& \nabla \varphi(x) = \sum_{i=1}^{m}\varphi'_{i}(a_{i}^{\T}x)a_{i},\ \ \nabla^{2}\varphi(x) = \sum_{i=1}^{m}\varphi''_{i}(a_{i}^{\T}x)a_{i}a_{i}^{\T}\\
	& \nabla^{2} \varphi(x)u = \sum_{i=1}^{m}\varphi_{i}''(a_{i}^{\T}x)a_{i}(a_{i}^{\T}u).
	\ena
	The computational complexities of $\nabla \varphi(x)$ and the HVP $\nabla^{2}\varphi(x)u$ are of the same as $O(nm)$, while the complexity of $\nabla^{2}\varphi(x)$ is $O(n^{2}m)$.\qed
\end{exam}

If we select $u$ in \eqref{equ:broyden_update_formula} greedily from the standard unit vectors of $\Rn$, i.e., 
\bea\label{equ:greedy_selected_u}
g_{A}(G)=  \underset{u\in\lbrace e_{1}, \ldots, e_{n} \rbrace}{\argmax} \frac{\langle Gu, u \rangle}{\langle Au, u \rangle},
\ena
 the  GQN update is iteratively given by
\bea\label{iter_matrix}
G^{k+1}& = \brd (G^{k}, A, g_{A}(G^{k}))\\
&\equdef  \gbd(G^{k}, A).
\ena
%\begin{rem}\label{gredyvec}
%In this work, we restrict the greedy vector to be a standard unit vector, leading to that the %computation of the matrix-vector product $Ag_{A}(G)$ in \textsf{gBroyd} only involves one %column of $A$.  If the greedy vector is arbitrarily  selected from any set of basis vectors of %$\Rn$ in  \citet{rodomanov2021greedy}, then $Ag_{A}(G)$ involves computing multiple %columns of $A$, which might violate the sprit of QN methods. 

% to save the computational cost.  Since $Ae_i$  exactly returns the $i$-th column of $A$, we must use some column of $A$ to compute the \textsf{gBroyd}. 
%\end{rem}
A striking feature of the GQN is that we are able to explicitly evaluate the convergence rate of $
\{G^{k}\}$ to $A$ via the following metric
\bea\label{equ:G_measure_defi}
\sigma_{A}(G)=\langle A^{-1}, G-A \rangle = \langle A^{-1}, G\rangle - n
\ena
where $n$ is the dimension of the vector space.
%Then,  $ G^{k}$ in \eqref{iter_matrix} will approach $A$ at a linear rate which is crucial to our algorithm design.
\begin{lem}[\citet{rodomanov2021greedy}]\label{lem:GQN_update_linear_rate}
If 	$
	0 \prec \mu I \preceq A \preceq \omega I
	$
and $G^k \succeq A $ in \eqref{iter_matrix},  then 
	\bea\label{equ:GQN_update_linear_rate}
	& \sigma_{A}(G^{k+1}) \leq (1 - \mu / (n\omega)) \sigma_{A}(G^{k}).
	\ena
\end{lem}

Under suitable conditions, the GQN method using a correction strategy has been shown to achieve a local explicit superlinear rate in \citet{rodomanov2021greedy}. In this work, we design an adaptive GQN method to achieve {\em global} convergence in Section \ref{sec1_global_convergence} and extend it to the master-worker setting in Section \ref{sec:GQN-LA}. Both can also achieve local linear and superlinear convergence. Moreover, the communication complexity of the distributed version is comparable to that of distributed first-order methods. 

\section{The adaptive greedy quasi-Newton method in the centralized setting}\label{sec1_global_convergence}
In this section, we first propose a novel adaptive GQN (AGQN) method where a simple mechanism is designed to adaptively adjust stepsizes to ensure its global convergence.  Our key idea lies in the use of Lemma \ref{lem:GQN_update_linear_rate} to adaptively tune stepsizes to uniformly bound the Hessian approximation error. Then, we establish explicit non-asymptotic bounds for both the global convergence rate and local superlinear rate of the AGQN. 

Consider an unconstrained optimization problem
\bea\label{equ:obj}
\underset{x\in\Rn}{\minimize}\ \  f(x).
\ena
Suppose that $x$ is the current iterate and $G\succ0$ is an estimate of $\nabla^2 f(x)$. A generic QN method has the following form
\bea\label{equ:x_update}
x^{+} = x - \alpha \cdot G^{-1}\nabla f(x)
\ena
where $0<\alpha\le 1$ is the stepsize and $f$ satisfies the following condition. 
\begin{assum}\label{assum:f}
	\begin{enumerate}[label=(\alph*)]
		\item $f$ is \textit{$\mu$-strongly convex and $\omega$-smooth, i.e., }
		\bea\label{equ:assum_1}
		\mu I \preceq \nabla^2 f(x) \preceq \omega I, \forall x\in\Rn
		\ena
		in which case we denote its condition number by $\kappa=\omega/\mu$.
		\item The Hessian of $f$ is \textit{$L$-Lipschitz continuous:}
		\bea\label{equ:assum_2}
		\Vert \nabla^2 f(y) - \nabla^2 f(x) \Vert \leq L\Vert y - x \Vert, \forall x,y\in\Rn.
		\ena
	\end{enumerate}
\end{assum}
Under Assumption \ref{assum:f}, one can establish the so-called strong self-concordancy of $f$.
\begin{lem}(\citet{rodomanov2021greedy})\label{lem:self_concordance}
	Suppose that $f(x)$ satisfies Assumption \ref{assum:f}. Then $f(x)$ is strongly self-concordant with a constant $M \leq L/\mu^{3/2} > 0$, i.e.,  for all $x,y,z,w \in \Rn$, it holds that          
	\bea\label{equ:aself_concordant}
	\nabla^2 f(y) - \nabla^2 f(x) \preceq M \Vert y - x \Vert_{
	\nabla^{2}f(z)} \nabla^2 f(w).
	\ena
\end{lem}
%\begin{pf}
%By $L$-Lipschitz continuity of the Hessian, we have that
%\bea
%&\hspace{-1cm} \nabla^{2}f(y) - \nabla^{2}f(x)
%\preceq L\Vert y- x \Vert I   \\
%& = L ( (y-x)^{\tp}(y-x))^{1/2} I \\
% & \preceq L ( (y-x)^{\tp} \nabla^{2}f(z) (y-x))^{1/2} I / \mu^{1/2}\\
% &= L \Vert y-x\Vert_{\nabla^{2}f(z)} I / \mu^{1/2}\\
%& \preceq L \Vert y-x\Vert_{\nabla^{2}f(z)} \nabla^{2}f(w) / \mu^{3/2}
%\ena 
%where the second and third inequalities follow from $\mu$-strongly convexity of  $\nabla^{2}f(z)$ and $\nabla^{2}f(w)$. Let $M = L / \mu^{3/2}$ and the proof is finished. 
%\end{pf}

	For convenience we shall directly use $M$ as the constant of strong self-concordancy in the sequel. Assumption \ref{assum:f} is common in the analysis of Newton-type methods \citep{polyak2020new, Boyd2004}. Note that the strong self-concordancy is the key for the correction strategy to \citet{rodomanov2021greedy}.

\subsection{Adaptive stepsizes for the global uniform boundedness of  Hessian estimate errors}\label{subsec:3_1}
%The key for the global convergence of QN methods is to choose stepsizes to ensure objective function descent per iteration, which is usually achieved by a series of backtracking line searches, e.g., finding a stepsize $\alpha=c \cdot \rho^k$ to satisfy two conditions \cite{byrd1987global}:
%\bea\label{equ:quasi_newton_line_search}
%& f(x + \alpha d) \leq f(x) + \rho \alpha \nabla f(x)^{\T}d\\
%&\nabla f(x + \alpha d)^{\T}d \geq \beta \nabla f(x)^{\T}d
%\ena
%where $d = - G^{-1}\nabla f(x)$ is a descent direction, $0 < \rho < 1/2$ and $\rho < \beta < 1$. However, in the distributed setting, such a line search is exceedingly time-consuming since it has to evaluate global objective function $f$ and its gradient multiple times per iteration. 
%
%An alternative is to utilize more elaborate properties of $f(x)$ such as Lipschitz Hessian and strongly convexity, to design an adaptive stepsize in a compact form. For example, Polyak \textit{et al.} \cite{polyak2020new} derive an adaptive stepsize for the exact Newton method and obtain global quadratic convergence. Inspired by this, we shall attempt to derive an adaptive mechanism to choose stepsizes for the GQN method.

The key idea for our AGQN is that the Hessian approximation error can be explicitly bounded in the optimization progress. 
Specifically, suppose that $G\succeq \nabla^2 f(x)$ is a ``good" approximation of $\nabla^{2}f(x)$, i.e., there exists a small constant $\varepsilon \geq 0$ such that
\bea\label{equ:G_condition}
\sigma_{\nabla^2 f(x)} (G) \leq \varepsilon.
\ena

To approximate $\nabla^{2}f(x^{+})$,  a corrected $G$ is employed in the GQN, i.e., 
\bee\label{correction}
\tG^+ = (1 + Mr^+)G, ~\textrm{and~} G^{+} =\gbd(\tG^+, \nabla^{2}f(x^{+}))  \\                   
\ene 
where  $\gbd$ is defined in \eqref{iter_matrix} and 
\bee
r^+= \Vert x^{+} -x \Vert_{\nabla^{2}f(x)}. \label{residual}
\ene
Note that $r^+$ can also be efficiently computed by the HVP (see Example \ref{example:1}). Such a correction strategy is essential as it not only ensures an upper Hessian estimate but also the error can be evaluated. Specifically, it follows from \citet[Lemma 4.8]{rodomanov2021greedy} that
 \bea\label{equ:error_G+}
&\sigma_{\nabla^{2}f(x^{+})}(G^+) \\
\quad & \leq \left(1 -\frac{1}{n\kappa}\right)(1 + Mr^+)^2\left(\varepsilon + \frac{2nMr^+ }{1 + Mr^+}\right).
\ena
Apparently, the Hessian approximation error at $x^+$ grows as $r^+$ increases. In the global analysis where $x$ can be arbitrarily far away from the optimal solution, we cannot bound $r^+$ and thus the Hessian approximation error is not restricted, which is the main challenge to establish global convergence for the GQN method. 

To resolve it, we propose to use \textit{multiple} \textsf{gBroyd} before \eqref{correction} to refine our estimate of $G$. Specifically, let $G^{(0)} = G$ and
 \bea
 \label{refinehessian}  
G^{(\tau)}&=\gbd(G^{(\tau-1)}, \nabla^{2}f(x)) \\
&=\ldots=\gbd^{(\tau)}(G^{(0)}, \nabla^{2}f(x))           
\ena
where $\tau$ is a nonnegative integer and $\gbd^{(\tau)}$ denotes the $\tau$ times of compositions of $\gbd(\cdot, \nabla^{2}f(x))$. 

Let $G=G^{(\tau)}$  in \eqref{correction}, i.e., 
\bee\label{correction2}
\tG^+ = (1 + Mr^+)G^{(\tau)}, ~\textrm{and}~G^{+} =\gbd(\tG^+, \nabla^{2}f(x^{+})).  \\                   
\ene 
%Inspired by the linear convergence of GQN update in Lemma \ref{lem:GQN_update_linear_rate}, 
%Suppose that we apply $\tau$ times GQN update to $G$ and obtain $\barG^{\tau}$, i.e., 
%\bea\label{equ:tG_defi_1}
%\barG^{l} = \brd_{\phi^{l}} (\barG^{l - 1}, \nabla^{2}f(x), \baru_{x}(\barG^{l - 1})), l=1,\cdots,\tau
%\ena
%where $\barG^{0} = G$. For simplicity, we directly write it as 
%\bea\label{equ:tG_defi_2}
%\barG^{\tau} = \brd^{(\tau)}(G, \nabla^{2}f(x), \baru_{x}(G)).
%\ena
Jointly with Lemma \ref{lem:GQN_update_linear_rate}, it holds that
\bea\label{equ:Broyd_inequ2}
	&\sigma_{\nabla^{2}f(x^{+})}(G^+)\leq \left(1 - {1}/{(n\kappa)}\right)\cdot  \\
	 &\hspace{1cm} \left(1 + Mr^+\right)^2\left(\left(1 - \frac{1}{n\kappa}\right)^\tau \varepsilon + \frac{2nMr^+}{1 + Mr^+}\right) 
	\ena
which is a tighter upper bound than \eqref{equ:error_G+}. 

Since $r^+$ in \eqref{equ:Broyd_inequ2} can be upper bounded by the stepsizes, we are able to adaptively tune them to achieve a uniform upper bound of  \eqref{equ:Broyd_inequ2} in the form of \eqref{equ:G_condition}. We formally state this result whose proof is given in Appendix \ref{app:proof_uniform_bound}. 

%, we have that the approximation error of $\barG^{\tau}$ decreases linearly, i.e., 
%\bea\label{equ:G_tau}
%\nabla^2 f(x) \preceq \bar{G}^{\tau} \preceq \left(1 + (1 - \frac{\mu}{n\omega})^{\tau}\varepsilon\right) \nabla^2 f(x).
%\ena
%Then we set $\tG = (1 + Mr)\barG^{\tau}$ and update $\tG$ into $G^{+}$ by Step 3 of Alg. \ref{alg:greedy_quasi_newton}. We hope that the decreased error of $\barG^{\tau}$ can compensate the increased error of $G^{+}$ in \eqref{equ:error_G+} if $r$ is sufficiently small,  which can be done by choosing a sufficiently small stepsize $\alpha$ and update $x$ by
%
%The following lemma verifies this intuition. It shows that under a sufficiently small stepsize, the error bound of $G^{+}$ will not be enlarged compared with $G$ after one iteration. The proof is in the appendix.
\begin{lem}\label{lem:uniform_bound}
	Consider the QN method of \eqref{equ:x_update} where $G$ and its GQN update satisfy \eqref{equ:G_condition} and  \eqref{correction2} respectively. Let 
	\bee\label{equ:bar_alpha_defi}
	\alpha_\tau(x) =\frac{ c_{\tau, \varepsilon} } {M \Vert \nabla f(x) \Vert_G^* } 
	\ene
	where $\tau$ is a nonnegative integer,
	\bea\label{equ:bar_r_tau}
	c_{\tau, \varepsilon}= - \frac{\sigma_\tau} {\sigma_\tau + n} + \frac{\sqrt{\rho \sigma_{\tau}^2 + (\sigma_\tau + n)(1 - \rho^{\tau + 1})\varepsilon}} {\sqrt{\rho} (\sigma_\tau + n)},
	\ena
	 $\rho = 1 - 1/ (n\kappa)$ and $\sigma_\tau = \rho^\tau \varepsilon + n$.
	If $f$ satisfies Assumption \ref{assum:f} and 
 $\alpha \le \alpha_\tau(x)$, then it holds that
	\bea\label{equ:G_error_inequ}
	&\sigma_{\nabla^2 f(x^+)}(G^+) \leq \varepsilon.
	\ena
\end{lem}
In view of \eqref{equ:G_condition} and \eqref{equ:G_error_inequ},  the Hessian approximation error is uniformly bounded, which is not covered by the local results in \citet{rodomanov2021greedy}. To capture the dependence of \eqref{equ:bar_r_tau} on the main constants, one can easily derive a lower bound as $c_{\tau, \varepsilon} \geq \Omega( (\tau + 1) \varepsilon / ( (n + \varepsilon) (\kappa n + \tau)))$, where $\Omega(\cdot)$ is Big Omega notation.

By \eqref{equ:error_G+} and \eqref{equ:Broyd_inequ2}, we use \eqref{refinehessian} to refine our estimate of $\nabla^{2}f(x)$ and adaptive stepsizes to achieve \eqref{equ:G_error_inequ}.  Intuitively, a larger $\tau$ will lead to a smaller Hessian approximation error (cf. \eqref{equ:Broyd_inequ2}) and possibly a larger objective reduction, though it also induces a higher computational cost. %Specifically, each \textsf{gBroyd} in \eqref{refinehessian} requires to compute a HVP which retains a column of $\nabla^2 f(x)$. 
	 On the other hand, if $\tau$ is small, it follows from \eqref{equ:bar_r_tau} that the allowable stepsize is small as well.  In this view, $\tau$ can be considered as a tradeoff factor among the  performance of approximating Hessian, the computational overhead and optimization stepsize. In practice we suggest to choose $\tau \ll n$.

\subsection{Adaptive stepsizes for successive improvements and the AGQN method}\label{subsec:3_2}
Till now, we have established the uniform boundedness of Hessian estimate errors, which is key to the design of adaptive stepsizes to achieve successive improvements of the objective function. 
As in \citet{polyak2020new}, we adopt the following metric to measure the algorithm performance
\bee\label{per_metric}
\beta(x) = L \Vert \nabla f(x)\Vert / (2\mu^{2}).
\ene

%	Note that the right side of \eqref{equ:beta_inequ} is a quadratic function on $
%	\alpha\in(0, \bar{\alpha}_{1}]$ where $\bar{\alpha}_{1} \equdef \min\lbrace \bar{\alpha}(\tau), 1\rbrace$. We minimize the right side with respect to $\alpha$:
%	\begin{subequations}\label{equ:beta_minimize}
%		\begin{align}
%			& \textrm{When}\ \ 1 / (2(1 + \varepsilon)\beta) \leq \bar{\alpha}_{1},\\
%			&\quad\quad\quad \alpha^{*} = 1 / (2(1 + \varepsilon)\beta),\ \beta^{+} \leq \beta - 1 / (4(1 + \varepsilon^{2})).  \\
%			& \textrm{Otherwise}, \alpha^{*} = \bar{\alpha}_{1},\ \beta^{+} \leq 1 - \bar{\alpha}_{1} / (2(1 + \varepsilon)))\beta.
%		\end{align}
%	\end{subequations}

% for $G$ in Lemma \ref{lem:uniform_bound}, we can develop an adaptive rule to choose stepsize and ensure objective improvment per iteration. We give it in the following lemma and the proof is in the appendix. Our technique is slightly different from \cite{polyak2020new} since we have to deal with approximate Hessian instead of the exact Hessian. 
\begin{lem}\label{lem:adp_step}  Suppose that all the conditions in Lemma \ref{lem:uniform_bound} with 
	\bea\label{equ:varepsilon_0}
	\varepsilon \leq \varepsilon_{0} \equdef 1 / (2 \kappa - 1)
	\ena
hold. 
Let 
\bee
\bar{\alpha}_\tau(x) = \min \{ \alpha_\tau(x) , 1\},
\ene
\bee\label{equ:alpha_beta_defi}
\alpha(x)= \min\{\alpha_\tau(x), 1 / (4 \beta(x)), 1 \}.
\ene
If $f$ satisfies Assumption \ref{assum:f}, then it holds that
%	Let $x\in \Rn$ and $G\succ 0$ satisfying \eqref{equ:G_condition}. $x^{+}$ is defined by \eqref{equ:x_update}. 
%	Define $$\beta \equdef L \Vert \nabla f(x)\Vert / (2\mu^{2}), \quad  \beta^{+} \equdef L\Vert \nabla f(x^{+})\Vert / (2\mu^{2}).$$
%	Also define 
%	\bea\label{equ:bar_alpha_1_defi}
%	\bar{\alpha}_{1} \equdef \min\lbrace 1, \bar{\alpha}(\tau)\rbrace, 	\baralp_{2} \equdef \frac{1}{2(1 + \varepsilon)\beta}, \bar{\beta}_{2}(\tau) \equdef \frac{1}{2(1 + \varepsilon)\bar{\alpha}_{1}}.
%	\ena 
%	where $\baralp(\tau)$ is defined in \eqref{equ:bar_alpha_defi}.
%	If we choose $\alpha$ as 
%	\bea\label{equ:adp_step_choice_1}
%	\alpha = \min\lbrace 1, \baralp(\tau), \baralp_{2} \rbrace, 
%	\ena
%	we can obtain the follwing objective improvement:
	\bea\label{equ:adp_step_descent}
	& \beta(x^+) \leq
	\begin{cases} 
		\beta(x)- 1 / 16, \quad\quad \text{if}~\beta(x) > 1 / (4 \bar{\alpha}_\tau(x) ), \\
		\
		(1 - \bar{\alpha}_\tau(x) / 4 )\beta(x), \quad \text{otherwise}.
	\end{cases}
	\ena
\end{lem}
The proof is given in Appendix \ref{app:proof_adp_step}. By telescoping the definition of $\alpha_{\tau}(x)$ in \eqref{equ:bar_alpha_defi} and using $G \succeq \nabla^2 f(x) \succeq \mu I$, we have that $\alpha_{\tau}(x) \geq \sqrt{\mu}c_{\tau, \varepsilon_0} / (M\Vert \nabla f(x)\Vert)$, which is increasing as $\beta (x) $ decreases. This means that $\bar{\alpha}_{\tau}(x)$ has a monotonically increasing lower bound until it reaches 1 and thus $\beta(x)$ can be successively reduced.

To sum up, our adaptive stepsizes are explicitly given by
\bea\label{equ:adap_stepsize}
\alpha(x) = 
\min\{ \alpha_{\tau}(x), 1/(4\beta(x)), 1 \}
\ena
where $\alpha_{\tau}(x)$ and $\beta(x)$ are defined in Lemmas \ref{lem:uniform_bound} and \ref{lem:adp_step}, respectively. Then, we obtain a novel Adaptive Greedy Quasi-Newton (AGQN) method in Alg. \ref{alg:AGQN}. Jointly with Lemma \ref{lem:uniform_bound} and an initial condition $\sigma_{\nabla^2f(x^0)}(G^0) \leq \varepsilon_{0}$, Lemma 5 actually establishes global convergence of the AGQN since $\sigma_{\nabla^2f(x^k)}(G^k) \leq \varepsilon_{0}$ holds for all $k$ and $\beta(x^k)$ will be decreased successively. 

\begin{algorithm}[!t]
	\caption{The line search free AGQN} \label{alg:AGQN}
	\begin{itemize}[leftmargin=*, label=$\bullet$]
		\item{\bf Initialization:} Choose $x^{0}$, $G^{0}$ and $\tau\geq 0$.
		\item{\bf For $k \geq 0$ iterate:}
		\begin{enumerate}[leftmargin=*, ]
			\renewcommand{\labelenumi}{\theenumi.}
			\renewcommand{\labelenumii}{(\arabic{enumii})}
			\item Compute $\nabla f(x^{k})$ and $\beta(x^{k})$ using \eqref{per_metric}. 
			\item Select an adaptive stepsize $\alpha^{k}=\alpha(x^k)$ by \eqref{equ:adap_stepsize} and update	
			$$x^{k+1} = x^{k} - \alpha^{k}(G^{k})^{-1}\nabla f(x^{k}).$$ 
			\item
			\begin{enumerate}[label=(\alph*)]
				\item Perform $\tau$ times of \textsf{gBroyd} updates on $G^{k}$ to refine our estimate of $\nabla^{2}f(x^{k})$, i.e., 
				$$\underline{G}^{k}=\textrm{gBroyd}^{(\tau)}(G^{k}, \nabla^{2}f(x^{k})).$$ 				
				\item Perform a \textsf{gBroyd} update to obtain an upper estimate of $\nabla^{2}f(x^{k+1})$, i.e., 
				$$G^{k+1} = \textrm{gBroyd}(\overline{G}^{k+1}, \nabla^{2}f(x^{k+1}))$$
				where 
				$\overline{G}^{k+1} = (1 + Mr^{k+1})\underline{G}^{k}$ and $r^{k+1} = \Vert x^{k+1} - x^{k} \Vert_{\nabla^2 f(x^{k})}.$
			\end{enumerate}
		\end{enumerate}
		\item{\bf Until} a termination condition is satisfied.
	\end{itemize}
\end{algorithm}

Note that in Alg. \ref{alg:AGQN}, all the terms involving Hessian can be computed via HVPs. Besides the use of adaptive stepsizes that is different from the GQN in \citet{rodomanov2021greedy}, we further perform $\tau$ updates of \textsf{gBroyd} for $G$ to refine our estimate of $\nabla^{2}f(x)$ for the global convergence.  Since $\tau$ can be as small as one (see Lemma \ref{lem:uniform_bound}) and the cost of computing the HVP $\nabla^{2}f(x)g$ is similar with $\nabla f(x)$ (see Example \ref{example:1}), this extra computation cost can be acceptable. Also note that while \textsf{gBroyd} is not given in the form of directly updating $(G^k)^{-1}$, one can still compute $(G^k)^{-1}$ with an $O(n^2)$ complexity by exploiting the low-rank structure of \textsf{gBroyd} and invoking the Sherman-Morrison-Woodbury formula \citep{horn2012matrix}.

\begin{rem} 
	The AGQN involves several constants $\mu$, $\omega$ and $L$ that are commonly used in the theoretical analysis of second-order methods \citep{rodomanov2022rates,jin2022non,rodomanov2021new}. It is possible to design algorithms to estimate these parameters. For example,  \citet{polyak2020new} proposes an adaptive Newton method to deal with the case of unknown $\mu$, but the cost is that the Newton direction is evaluated multiple times per iteration.  There are also methods of estimating the smoothness parameters in \citet{mishchenko2021regularized} and \citet{Malitsky2020}. In practice, $f$ may have some  specific structure, e.g., the logistic regression problem in Section \ref{sec:experiments}, one can easily derive their bounds and then resort to manual tuning to improve performance if needed. Though such a problem is nontrivial, we are short of space for further investigation and use them directly to present our results. 
\end{rem}

%When the iterate comes close to the optimal solution, the adaptive stepsize $\alpha$ is simply set to one, after which we use the GQN in \citet{rodomanov2021greedy} to give the local superlinear convergence result. 

%To elaborate it, we use \citet[Theorem 4.9]{rodomanov2021greedy} to conclude that the superlinear convergence occurs if
%\bea\label{equ:superlinear_neighbor}
%\lambda(x) \leq \frac{\ln(2)}{4(2n + 1)M\kappa}
%\ena
%where $\lambda(x)$ is given in \eqref{lambdaf}. Since $\lambda(x)$ depends on $\nabla^{2}f(x)$ and thus is not available, we use an upper bound via the $\mu$-strongly convexity in Assumption \ref{assum:f}(b) to obtain that
%\bea
%\lambda(x) \leq \Vert \nabla f(x)\Vert / \mu^{1/2}.
%\ena
%In light of \eqref{per_metric} and \eqref{equ:superlinear_neighbor}, a sufficient condition for the local superlinear convergence is explicitly given as 
%\bee\label{equ:beta_s_defi}
%\beta(x) \leq \frac{\ln(2)L}{8(2n + 1)\mu^{1/2}M\omega}\equdef \beta_s.
%\ene
%Since this condition is conservative,  
%the superlinear rate can possibly occur even if it is not satisfied. 

\subsection{Local convergence rates of the AGQN}\label{subsec:3_3}
When the iterate of the AGQN comes close enough to the optimal solution, the adaptive stepsize in \eqref{equ:adap_stepsize} simply becomes one, after which we can discuss the local convergence behaviour of the AGQN. We define the following suboptimality metric for $f$
\bee\label{lambdaf}
\lambda(x) =  ( \nabla f(x)^{\tp} \nabla^{2}f(x)^{-1}\nabla f(x) )^{1/2}.
\ene
Based on this, we give the following theorem to establish local convergence properties for the AGQN. It shows that, with bounded Hessian approximation error per iteration, the AGQN achieves a local constant linear rate and a local superlinear rate. The proof is in the Appendix \ref{app:AGQN_local}.

\begin{thm}\label{thm:AGQN_local}{\bf(Local convergence of AGQN)} Suppose that Assumption \ref{assum:f} holds and $\sigma_{\nabla^2f(x^0)}(G^0) \leq \varepsilon_{0}$ (cf. \eqref{equ:varepsilon_0}). Suppose that $x^0$ is sufficiently close to the solution, i.e.,
		\bea\label{equ:AGQN_ll_cond1}
		\Vert \nabla f(x^0) \Vert \leq \min \lbrace \dfrac { \sqrt{\mu}c_{\tau, \varepsilon_{0}} } { M }, \dfrac{\mu^2}{2L}  \rbrace.
		\ena
		If 
		\bea\label{equ:AGQN_ll_cond2}
		M\lambda(x^0)\leq 1/4,
		\ena
		we obtain the following linear rate
		\bea\label{equ:AGQN_ll_rate}
		\lambda(x^{k+1}) \leq \dfrac{3}{4}\lambda(x^k).
		\ena
		Furthermore, if 
		\bea\label{equ:AGQN_ls_cond2}
		M\lambda(x^0)\leq \dfrac{\ln 2} { 8(6n + 1) },
		\ena
		 we obtain the following superlinear rate
		\bea\label{equ:AGQN_ls_rate}
		\lambda(x^{k+1}) \leq \left( \max \lbrace \rho^{\tau + 1}, \dfrac{3}{4}\rho \rbrace \right)^k 2\bar{\Phi}_0 \cdot \lambda(x^k)
		\ena
		where $\rho$ is in Lemma \ref{lem:uniform_bound} and $ \bar{\Phi}_0= \varepsilon_0 + n\ln 2 / (12n  + 2 )$.
		
\end{thm}

It is worth noting that both rates are better than those in \citet{rodomanov2021greedy}, which is natural since we conduct $\tau$ additional \textsf{gBroyd} updates per iteration. Moreover, the superlinear rate in \eqref{equ:AGQN_ls_rate} can get better as $\tau$ increases, and if $\tau = 0$ it also reduces to essentially the same rate in Theorem 4.9 in \citet{rodomanov2021greedy}. Combining the global convergence guarantee in Lemma \ref{lem:adp_step} and local convergence results in Theorem \ref{thm:AGQN_local}, we can establish the total iteration complexity of AGQN in the following corollary. Its proof is given in Appendix \ref{app:cor_AGQN}.
\begin{cor}\label{cor:AGQN}
	Suppose that Assumption \ref{assum:f} holds and $\sigma_{\nabla^2f(x^0)}(G^0) \leq \varepsilon_{0}$. The totol iteration complexity for the AGQN in \ref{alg:AGQN} to achieve $\Vert \nabla f(x) \Vert \leq \epsilon$ is given by
	\bea\label{equ:AGQN_comp}
	 & \widetilde{\bigO} \left( \frac {nL (\kappa n + \tau) \kappa }{ \mu^2 \tau}  \Vert \nabla f(x^0) \Vert \right) \\
	 & \quad\quad\quad\quad + \bigO(\sqrt{\min \lbrace 4n\kappa/3, n\kappa + \tau \rbrace } \sqrt{\log(\bar{\Phi}_0 / \epsilon)}). 
	 \ena
\end{cor}
Different from FedNL \citep{Safaryan2022} and GIANT \citep{Wang2018}, the superlinear rate in Theorem \ref{thm:AGQN_local} depends on the condition number $\kappa$ and the dimension of decision vector $n$.  This is perhaps due to that AGQN only has access to an {\em incomplete} Hessian via a limited number of HVPs while FedNL computes the full Hessian before compression and the GIANT solves a subproblem involving with the complete Hessian. Kindly note that even for the standard quasi-Newton methods, e.g., SR1, DFP and BFGS, the non-asymptotic {\em local} convergence rate was unknown until very recently, which also depends on the condition number, see eq. (50) in \citet{rodomanov2021new} and eq. (73) in \citet{jin2022non}. Indeed, it is interesting to design methods to overcome such a limitation.

\section{The distributed AGQN in the master-worker setting } \label{sec:GQN-LA}
In this section, we extend the AGQN to the master-worker setting, the purpose of which is to collaboratively solve the following distributed optimization problem
\bea\label{equ:intro_obj}
\underset{x\in\Rn}{\minimize}\ \  f(x)=\sum_{i=1}^{p}f_{i}(x)
\ena
where each worker $i$ has access to a local objective function $f_{i}$ and the master computes the optimal decision vector using information from workers. We suppose that Assumption \ref{assum:f} holds for each $f_i$ with $\mu$, $\omega$ and $L$. This is reasonable since if Assumption \ref{assum:f} holds for $f_i$ with $\mu_i$, $\omega_i$ and $L_i$, we can simply take $\mu = \min_i \mu_i$, $\omega_i = \max_i \omega_i$ and $L = \max_i L_i$. Then the following lemma summarizes properties of $f$ and its proof is given in the Appendix \ref{app:f_prop}.
\begin{lem}\label{lem:f_prop}
	Suppose that Assumption \ref{assum:f} holds for each $f_i$ with $\mu$, $\omega$ and $L$. Let $M$ be the constant of strong self-concordance of each $f_i$. Then we have the following properties of $f$:  
	\begin{enumerate}[label=(\alph*)]
		\item $f$ is $p\mu$-strongly convex and $p\omega$-smooth.

		\item The Hessian of $f$ is $pL$-Lipschitz continuous.
	
	    \item $f$ is strongly self-concordant with constant $M_f$ where $M_f \leq \sqrt{p} M$.
	\end{enumerate}
\end{lem}

\subsection{The bottleneck of implementing the AGQN in the master-worker setting}\label{subsec:4_1}
It is nontrivial to implement the AGQN in the master-worker setting.  
For the master to complete one \textsf{gBroyd} in Step 3(a), it follows from \eqref{iter_matrix} that  the necessary information includes a greedy vector $g:=g_{\nabla^2 f(x)}(G)$ and a HVP  $\nabla^2 f(x) g$. By \eqref{equ:greedy_selected_u}, the computation of $g$ and $\nabla^2 f(x) g$ requires the master to firstly pull all the diagonal elements of $\nabla^2 f_i(x)$ and then  $\nabla^2 f_i(x)g$ from every worker $i$.  Since $g$ is computed in the master, the worker $i$ cannot send them in a single packet. 
Thus, one update of \textsf{gBroyd} inevitably requires two rounds of coordinated communication between each worker and the master. 
By Step 3(a) and (b) in Alg. \ref{alg:AGQN}, {\em distributedly} updating $G$ requires to sequentially compute $(\tau + 1)$ updates of \textsf{gBroyd}  over the network. This implies that the number of communication rounds may become prohibitive in practice.    

\subsection{The distributed update of \textsf{gBroyd}}\label{subsec:dis_GQN}

In this subsection, we propose a ``divide-and-combine" strategy to ensure that each worker $i$ only needs to firstly pull the latest iterate from the master and then push a local vector $\cI_i$ of size $O((\tau+1) n)$ to the master. Since $\tau$ can be as small as one, the size of $\cI_i$ is an absolute constant number of times of that in the distributed first-order methods.  However, our distributed algorithm still retains the superb convergence property of the centralized AGQN.

Our central idea of ``divide" is that each worker $i$ maintains a local Hessian estimate $G_i$ of $\nabla^2 f_i(x)$ and extracts an $f_i$-dependent local vector $\cI_i$ for updating $G_i$ via \textsf{gBroyd}s. Then, the worker pushes $\cI_i$ to the master to help it to reconstruct the same copy of $G_i$ as that of the worker $i$\footnote{Though this potentially results in a scalability issue, it can be resolved (cf. Remark \ref{scalar}).}, based on which the master  simply ``combines" these $G_i$s to  form an estimate of   $\nabla^2 f(x)$ via $G=\sum_{i=1}^p G_i$. There are two remaining problems in such a natural idea, including (a) how to extract the local vector  $\cI_i$ and (b)  the major result of Lemma \ref{lem:GQN_update_linear_rate} is no longer available as the iterate of the form \eqref{iter_matrix} is destroyed under the ``divide-and-combine" policy.  We shall address the first problem in this subsection and leave the second to Section \ref{sec:converge}. 

Once the worker $i$ pulls the latest iterate $x^+$ from the master, it firstly refines its previous estimate of \ $\nabla^2 f_i(x)$ via $\tau$ updates of \textsf{gBroyd}, i.e., 
\bee
\underline{G}_i =  \gbd^{(\tau)}(G_i, \nabla^2 f_i(x)). \label{taugreedy}
\ene
By \eqref{equ:greedy_selected_u}, the computation of \eqref{taugreedy} requires to use $\tau$ vectors of HVPs with respect to $\nabla^2 f_i(x)$. Precisely, let the greedy vector be $g_i^{(t)}$ in the $t$-th update of \textsf{gBroyd} in \eqref{taugreedy}, then  $\nabla^2 f_i(x)g_i^{(t)}$ is the corresponding HVP.  Since $g_i^{(t)}$ is a standard unit vector of $\R^n$, it can be encoded by the position of its nonzero element and is denoted by an integer id$(g_i^{(t)})$.  Hence, the following information that is generated in the worker $i$ is sufficient for implementing \eqref{taugreedy}
\bee\label{infors}
\mathcal{S}_i=\left\{\begin{array}{ll}\emptyset &\text{if}~\tau =0,\\
\left\{\text{id}(g_i^{(t)}),\nabla^2 f_i(x)g_i^{(t)}\right\}_{t=1}^\tau &\text{if}~\tau>1.\end{array}\right.
\ene
That is, the worker $i$ uses $\underline{G}_i $ to refine its estimate of the Hessian matrix $\nabla^2 f_i(x)$, which is similar to Step 3(a) of the AGQN of Alg. \ref{alg:AGQN}.  

Then, it simply implements one \textsf{gBroyd} to estimate $\nabla^2 f_i(x^+)$, i.e., 
\bee
G_i^+=\gbd(\tG_i^+, \nabla^2 f_i(x^+)),~\text{and}~ \tG_i^+=(1 + Mr_i^+)\underline{G}_i\label{timeupdate}
 \ene
where $r_i^+=\|x^+ - x \|_{\nabla^2 f_i(x)}$. Similarly,  the \textsf{gBroyd} in \eqref{timeupdate} requires to use the index of a unit greedy vector $g_{i}$ and the associated HVP with respect to $\nabla^2 f_i(x^+)$. 

Overall, the information that is required for updating $G_i$ to $G_i^+$ via \eqref{taugreedy} and  \eqref{timeupdate} can be collectively summarized by a vector of size $\mathcal{O}((\tau + 1)n)$, i.e.,
\bea\label{equ:vec_pac_1}
 & \cI_i^+ = 
	  \lbrace\mathcal{S}_i,\nabla^{2}f_{i}(x^+)g_{i}^+, \textrm{id}(g_{i}^+), r_{i}^+ \rbrace,
\ena
which can be generated via a function $\textsf{Infvec}$ in Alg. \ref{alg:Gupdate} by the worker $i$. Since $\cI_i^+$ will be communicated to the master in a single packet, it can also update $G_i$ via  $\textsf{Gupdate}$ in Alg. \ref{alg:Gupdate}. Clearly, such a communication scheme is the same as that of distributed first-order methods, e.g., \citet{recht2011hogwild, lian2015asynchronous, leblond2017asaga, liu2015an, Du2021}. 
%The first part contains all the greedy vectors and HVPs that are required to update $G_i^-$ into $\underline{G}_i^{-}$ by $\tau$ sequential GQN updates with respect to $\nabla^2 f_i(x^-)$ (see step 3(a) in Alg. \ref{alg:AGQN}). The second line contains the greedy vector and HVP required to update $\underline{G}_i^{-}$ into $G_i$ by one GQN update with respect to $\nabla^2 f_i(x)$ (see step 3(b) in Alg. \ref{alg:AGQN}). After entering the superlinear phase, only the second part is required and we set $\tau = 0$ in $\cI_i$.

%With $\cI_i$, we define a function $\textsf{Gupdate}(\cdot)$ in Alg. \ref{alg:Gupdate} to describe how $G_i^-$ is updated into $G_i$. Note that all the vectors in $\cI_i$ are used by \textsf{gBroyd} and not explicited given.
\begin{algorithm}[t!]
	\caption{Computation of local information vector and Hessian estimate - from the view of worker $i$.} 
	\label{alg:Gupdate}
	\begin{itemize}[leftmargin=*, label=$\bullet$]
		\item {\bf function} $\cI_i^+=\textsf{Infvec}(G_i, x, x^+, f_i)$	 
		\begin{enumerate}[leftmargin=*, label=(\alph*)]
			\item Use \eqref{taugreedy} and \eqref{infors} to obtain 
			$\underline{G}_{i}$ and $\mathcal{S}_i$. 
			\item Perform \eqref{timeupdate} and \eqref{equ:vec_pac_1} to obtain $\cI_i^+$.
		\end{enumerate}
	\item {\bf function} $G_i^+=\textsf{Gupdate}(G_i, \cI_i^+)$
	\begin{enumerate}[leftmargin=*, label=(\alph*)]
			\item Use \eqref{taugreedy} and   \eqref{timeupdate}  to obtain $G_i^+$.
		\end{enumerate}
	\end{itemize}
\end{algorithm}

 \subsection{The distributed adaptive greedy quasi-Newton method}\label{subsec:4_3}
For the master-worker setting,  we propose the Distributed Adaptive Greedy Quasi-Newton (DAGQN) method in Alg. \ref{alg:DAGQN}, which utilizes the adaptive stepsize given in \eqref{equ:adap_stepsize_DAGQN} and the distributed update of \textsf{gBroyd} in Section \ref{subsec:dis_GQN}. Similar with AGQN in Alg. \ref{alg:AGQN}, the main computation cost for the master comes from $d^k = (G^k)^{-1} \nabla f(x^k)$. If using the technique of Remark \ref{scalar} and also invoking the Sherman-Morrison-Woodbury, one can compute $d^k$ with an $O(\min \lbrace p, n\rbrace n^2)$ complexity, which can be better than direct matrix inversion when $p < n$.
 
Since the \textsf{gBroyd} formula of \eqref{iter_matrix} is nonlinear, the relationship between $G^{k+1}$ and  $G^{k}$ in Alg. \ref{alg:AGQN} does not hold, meaning that we cannot use the key Lemma \ref{lem:GQN_update_linear_rate} to study the convergence of the DAGQN. Fortunately, we are still able to  establish explicit non-asymptotic bounds for the DAGQN by exploring the structure of Alg. \ref{alg:DAGQN}.

\begin{table*}[!htbp]
	\centering \caption{Comparison with existing distributed second-order methods.} 
	\begin{threeparttable}
	{\scriptsize 
	\begin{tabular}{|c|c|c|c|c|c|}
		\hline
		Methods & Types & Iteration complexity & \makecell{Comm. \\per round} & \makecell{Comp. per iterate \\ in a $\textrm{worker}^{(a)}$} & \makecell{Storage cost \\ in the master} \\ 
		\hline
		\makecell{GIANT}  & \makecell{local linear\\(GLM)} & $\bigO\left(\frac{\log(n\kappa / \epsilon)}{\log(m/n)}\right)$ & $\bigO(n)$ & Hessian of $f_i$& $\bigO(pn)$ \\  
		\hline
		$\textrm{NL}^{(b)}$  & \makecell{local superlinear \\ (GLM)} & $\bigO \left( \sqrt{\alpha_{\textrm{c}}} \sqrt{\log \left(LR\Phi^0 / \epsilon \right) } \right)$ & $ \bigO(n+*)  $ &  Hessian of $f_i$ & \makecell{$\bigO(n^2+ nmp)$}  \\
		%\hline
%		\makecell{DAN-LA \\ \citep{zhang2022distributed}} & \makecell{Global linear, \\ local spl, \\ implicit} & * & \makecell{ $\bigO(n)$ /\\ $\bigO(n + d_\mathcal{G})$ } & $\bigO(mn^2)$ & $\bigO(n^2)$ & \makecell{Decentralized, \\ compute $\nabla^2f_i(x)$ \\ explicitly}\\
		\hline
		\makecell{DANE}  & \makecell{global linear \\(for quadratics) } & $\bigO\left(\frac{\kappa^2}{m}\log\left(np/\delta\right)\log(1/\epsilon)\right)$ & $\bigO(n)$ & \makecell{solve an $f_i$-dependent  \\ subproblem}& $\bigO(n)$ \\
		\hline
		\makecell{DAGQN } & \makecell{global convergence, \\ local superlinear} & \makecell{$\bigO(\sqrt{\min \lbrace 4n\kappa/3, n\kappa + \tau \rbrace } \sqrt{\log(\bar{\Phi}_0 / \epsilon)})$} & $\bigO(\tau n)$ & \makecell{$\tau+1$ HPVs of $f_i$} & $\bigO(pn^2)$\\
		\hline 
	\end{tabular}}
\begin{tablenotes}\tiny
        \item [(a)] Except DANE, the cost of computing the gradient of $f_i$ is not included. 
	\item [(b)]  $\alpha_\textrm{c}$ depends on the compressor of NL,  $\Phi^0$ is the initial value of the potential function, $R=\max_{i,j}\Vert \xi_{ij}\Vert$ is the maximum norm of the training sample and $*$ denotes an implicit positive value. See \citet{Islamov2021} for details.  
\end{tablenotes}
	\end{threeparttable}
	\label{tab:methods_comp}
\end{table*}

\begin{savenotes}
	\begin{algorithm}[t!]
		\caption{The DAGQN: the global view} \label{alg:DAGQN}
		\begin{itemize}[leftmargin=*, label=$\bullet$]
			\item{\bf Initialization:}  Choose $x^{0}$ and $G_{i}^{0}, \forall i\in\{1,\ldots,p\}$.
			\item{\bf For $k \geq 0$ iterate:}
			
		 {\bf Each worker $i$ performs the following steps  in parallel.} 
				\begin{enumerate}[leftmargin=*, label=(\alph*)]
				\item If $k=0$, push $\nabla f_i(x^0)$ to the master.
				
				\item {\bf Pull communication:}~ Pull the latest iterate $x^{k+1}$ from the master.	
							
				\item  {\bf Local computation:} Invoke the \textsf{Infvec} function to update the local vector 
				$$\cI_i^{k+1}= \textsf{Infvec}(G_i^{k}, x^{k}, x^{k+1}, f_i) $$ 
				and the \textsf{Gupdate} function to update the local Hessian estimate 
				$$G_i^{k+1} = \textsf{Gupdate}(G_i^{k}, \cI_i^{k+1} ).$$ 
				\item {\bf Push communication:}~ Push the local information vector $\cI_i^{k+1}$ and gradient $\nabla f_i( x^{k+1} )$ to the master.
				\end{enumerate}
				
{\bf The master performs the following steps.} 
	\begin{enumerate}[leftmargin=*, label=(\alph*)]
		\item  {\bf Update the iterate:} Compute the update direction $d^k = (G^k)^{-1} \nabla f(x^k)$ . Use \eqref{equ:adap_stepsize_DAGQN} to determine an adaptive stepsize $\alpha^{k} = \alpha(x^{k})$. Then the decision vector is updated by 
				$$x^{k+1} = x^{k} - \alpha^{k} d^k.$$
				
		\item  {\bf Update Hessian estimate:} Invoke the \textsf{Gupdate} function to compute $G_{i}^{k+1} = \textsf{Gupdate}(G_{i}^{k}, \cI_i^{k+1})$. Aggregate the local Hessian estimates and gradients 
		$$G^{k+1}= \sum_{i=1}^{p} G_{i}^{k+1}, ~~~~\nabla f(x^{k+1}) = \sum_{i=1}^{p} \nabla f_{i}(x^{k+1}).$$
		\end{enumerate}
			\item{\bf Until} a termination condition is satisfied.
		\end{itemize}
	\end{algorithm}
\end{savenotes}
\section{Convergence analysis of the DAGQN}\label{sec:converge}
In this section, we study the convergence of the DAGQN and compare with the distributed second-order methods in the literature.  
\subsection{Explicit non-asymptotic bounds for the DAGQN}
Though  Lemma \ref{lem:GQN_update_linear_rate} cannot be used to study the convergence of DAGQN in Alg. \ref{alg:DAGQN},  this section proves that we can still give a similar convergence property of the  AGQN  and establishes explicit non-asymptotic bounds for its global convergence and local superlinear rate. 

First, we exploit the structure of \eqref{equ:intro_obj} and the distributed nature of DAGQN to propose a new metric to study how the Hessian approximation error grows. Specifically, let $G=\sum_{i=1}^pG_i$, $A=\sum_{i=1}^pA_i$ and $G_i\succeq A_i\succ 0, \forall i\in\mathcal{P}=\{1,\ldots,p\}$. The metric to quantify the gap between $G$ and $A$ is defined as 
\bee\label{equ:Delta:defi}
\Delta_{A}(G) \equdef \sum_{i=1}^{p} \sigma_{A_i} (G_i) = \sum_{i=1}^{p}\left(\langle  A_i^{-1}, G_{i}\rangle - n\right) 
\ene
under which Lemma \ref{lem:GQN_update_linear_rate} is also revised accordingly. 
\begin{lem}\label{lem:new_GQN_update_linear_rate}
Suppose that $0 \prec \mu I \preceq A_{i} \preceq \omega I$ and $G_i \succeq A_i $ for each $i\in\mathcal{P}$. Let $G^{+} = \sum_{i=1}^{p}G^{+}_{i}$ and 
  $G^{+}_{i} = \gbd(G_{i}, A_{i})$. Then, it holds that
	\bee\label{equ:G_global_bounds}
	 \Delta_{A}(G^{+}) \leq (1 - 1 / (n\kappa) ) \Delta_{A}(G).
	\ene
	where $\kappa = \omega / \mu$.
\end{lem}
\begin{pf} The proof can be easily derived by applying Lemma \ref{lem:GQN_update_linear_rate} to each $G_{i}, \forall i\in\mathcal{P}$ and then aggregating them. 
\end{pf}

In fact, $\Delta_{A}(G^k)$ is an upper bound of $\sigma_{A}(G^k)$ in \eqref{equ:G_measure_defi} (see Lemma \ref{lem:two_sigma_comparison} in Appendix \ref{app:proof_DAGQN_theo}).
Similar with Lemma \ref{lem:uniform_bound}, the following lemma shows that Hessian estimate error of the DAGQN can also be restricted under certain level in terms of the metric $\Delta_{\nabla^2 f(x)}(G)$.

\begin{lem}\label{lem:uniform_bound_Delta}
	For each $i\in\mathcal{P}$, suppose that Assumption \ref{assum:f} holds for $f_i$ with $\mu$, $\omega$ and $L$ and $G_i \succeq \nabla ^2 f_i (x)$. Suppose that $\Delta_{\nabla^2 f(x)}(G) \leq \varepsilon$. Let $x^+$ be defined by \eqref{equ:x_update}, and let $G^+ \equdef \sum_{i=1}^p G_i^+$ where $G_i^+$ is given by \eqref{timeupdate} for all $i\in\mathcal{P}$. Let 
	\bee\label{equ:bar_alpha_defi_Delta}
	\alpha_\tau(x) =\frac{ c_{\tau, \varepsilon} } {M \Vert \nabla f(x) \Vert_G^* } 
	\ene
	where $\tau$ is a nonnegative integer,
	\bea\label{equ:c_tau_Delta}
	c_{\tau, \varepsilon}= - \frac{\sigma_\tau} {\sigma_\tau + n\sqrt{p}} + \frac{\sqrt{\rho \sigma_{\tau}^2 + (\sigma_\tau + n\sqrt{p})(1 - \rho^{\tau + 1})\varepsilon}} {\sqrt{\rho} (\sigma_\tau + n\sqrt{p})},
	\ena
	$\rho = 1 - 1/ (n\kappa)$ and $\sigma_\tau = \rho^\tau \varepsilon + n\sqrt{p}$. 
	If $\alpha \le \alpha_\tau(x)$, then it holds that
	\bea\label{equ:G_error_inequ_Delta}
	&\Delta_{\nabla^2 f(x^+)}(G^+) \leq \varepsilon.
	\ena
\end{lem} 

Since $f$ in \eqref{equ:intro_obj} is $p\mu$-strongly convex and $pL$-Hessian Lipschitz, similar with \eqref{per_metric}, we define the following metric to measure the algorithm performance
\bee\label{per_metric_DAGQN}
\beta(x) = L \Vert \nabla f(x)\Vert / (2p\mu^{2}).
\ene
Then the adaptive stepsize for the DAGQN is given by
\bea\label{equ:adap_stepsize_DAGQN}
\alpha(x) = \min\{ \alpha_{\tau}(x), 1/(4\beta(x)), 1 \}
\ena
where $\alpha_\tau(x)$ is given in \eqref{equ:bar_alpha_defi_Delta} and $\beta(x)$ is defined in \eqref{per_metric}.
Similar with AGQN, we can establish non-asymptotic bounds for the DAGQN, containing global convergence guarantee, a local linear rate and a local superlinear rate. Its proof is deferred to the Appendix \ref{app:proof_DAGQN_theo}.

\begin{thm}\label{thm:DAGQN_theo} {\bf(Explicit non-asymptotic bounds of DAGQN)} For each $i\in\mathcal{P}$, suppose that Assumption \ref{assum:f} holds for $f_i$ with $\mu$, $\omega$ and $L$, and $G_i^0 \succeq \nabla ^2 f_i (x^0)$. Suppose that $\sigma_{\nabla^2f(x^0)}(G^0) \leq \varepsilon_{0}$ (cf. \eqref{equ:varepsilon_0}).

\begin{enumerate}[label=(\alph*)] 
\item {\bf (Global convergence)} The DAGQN in Alg. \ref{alg:DAGQN} converges at the following explicit non-asymptotic rates:
	\bea
	\hspace{-0.7cm}
	\beta(x^{k+1}) \leq  
	\begin{cases}
		\beta(x^{k}) - 1 / 16, ~ \text{if}~\beta(x^k) > 1 / (4 \bar{\alpha}_\tau(x^k)), \\
		(1 - \bar{\alpha}_\tau(x^k) / 4 )\beta(x^k), ~ \text{otherwise}
	\end{cases}
	\ena
	where $\bar{\alpha}_{\tau}(\cdot) = \min \lbrace \alpha_{\tau}(\cdot) ,1 \rbrace$ and $\alpha_{\tau}(\cdot)$ is defined in Lemma \ref{lem:uniform_bound_Delta}.  
	
\item {\bf (Local convergence)} Suppose that $x^0$ is sufficiently close to the optimal solution, i.e.,
\bea\label{equ:DAGQN_ll_cond1}
\Vert \nabla f(x^0) \Vert \leq \min \lbrace \dfrac { \sqrt{ p \mu}c_{\tau, \varepsilon_{0}} } { M }, \dfrac{p \mu^2}{2L}  \rbrace,
\ena
where $c_{\tau, \varepsilon_{0}}$ is defined in \eqref{equ:c_tau_Delta}. If
\bea\label{equ:DAGQN_ll_cond2}
\sqrt{p}M\lambda(x^0)\leq 1/4,
\ena
 we obtain the following linear rate
\bea\label{equ:DAGQN_ll_rate}
\lambda(x^{k+1}) \leq \dfrac{3}{4}\lambda(x^k).
\ena 
Furthermore, if
 \bea\label{equ:DAGQN_ls_cond2}
 \sqrt{p} M \lambda(x^0)\leq \dfrac{\ln 2} { 8(6n + 1) },
 \ena
  we obtain the following superlinear rate
 \bea\label{equ:DAGQN_ls_rate}
 \lambda(x^{k+1}) \leq \left( \max \lbrace \rho^{\tau + 1}, \dfrac{3}{4}\rho \rbrace \right)^k 2\bar{\Phi}_0 \cdot \lambda(x^k)
 \ena
 where $\rho = 1 - 1 / (n\kappa)$ and $ \bar{\Phi}_0= \varepsilon_0 + n\ln2/(12n + 2)$.
\end{enumerate}
\end{thm}
With Theorem \ref{thm:DAGQN_theo}, we can establish the iteration complexity of DAGQN is bounded by
\bea\label{equ:DAGQN_comp}
& \widetilde{\bigO} \left( \frac {n \sqrt{p} L (\kappa n + \tau) \kappa }{ \mu^2 \tau}  \Vert \nabla f(x^0) \Vert \right) \\
& \quad\quad\quad\quad + \bigO(\sqrt{\min \lbrace 4n\kappa/3, n\kappa + \tau \rbrace } \sqrt{\log(\bar{\Phi}_0 / \epsilon)} ). 
\ena
Its proof is similar with Corollary \ref{cor:AGQN} and we omit it due to space limitation. It is worth noting that compared with AGQN, our DAGQN has the same complexity for local superlinear convergence, which is the main concern of distributed second-order optimization.

\subsection{Comparison with the existing distributed second-order methods} 
	In the literature, the distributed methods are  specially developed to solve ERPs that the local function $f_i$ in \eqref{equ:intro_obj} is given in the following form
	\bea\label{equ:f_i_empirical}
		f_i(x) = \sum_{j = 1}^{m}  f_{ij}(x, \xi_{ij}), \ \  i\in\mathcal{P},
		\ena
	where $x\in\Rn$, $\xi_{ij}$ denotes a training sample and  $m$ is the sample size in each worker. Particularly, the GIANT  \citep{Wang2018}  and NL \citep{Islamov2021} focus on the GLM, i.e.,  
		\bea\label{equ:f_i_GLM}
		f_{ij}(x, \xi_{ij}) = \phi_{ij}(\xi_{ij}^{\mathsf{T}}x),
		\ena
	and the convergence result of the DANE \citep{Shamir2014} only holds for the case that $f_{ij}$ is  quadratic with a fixed Hessian. 
	
	Then, we compare the DAGQN with above-mentioned methods in terms of iteration complexity to reach an accuracy that $\|x^k-x^*\|\le \epsilon $, the number of scalars transmitted per round from each worker to the master,  the computational cost per iterate in a worker and the number of scalars that are stored in the master. The details are given in Table \ref{tab:methods_comp}. Since the major computational costs of the second-order methods are on the (quasi-)Newton direction, we only focus on the comparison of these costs.  Note that our storage cost is still less than NL when $m\gg n$, and our DAGQN requires only one communication round between the master and workers per iteration, which is not the case for the DANE and GIANT. 
	
 	It should also be noted that our line search free DAGQN  achieves global convergence for the general finite sum problem \eqref{equ:intro_obj}, which clearly covers the GLM in \eqref{equ:f_i_empirical}.

\begin{rem} \label{scalar}
To simplify presentation, the DAGQN of the present form has a storage cost $\bigO(pn^2)$ for the master. 
	A possible solution is that the master only stores $G^-$ with a storage cost $\bigO(n^2)$. By exploiting the low-rank structure of $\gbd$, each worker sends all the vectors and scalars that are required to compute the difference, namely $G_i-G_i^-$, to the master, which again takes a  storage cost $\bigO(n^2)$. Then, the master only needs to {\em sequentially} add these differences to $G^-$, leading to that $G=G^-+\sum_{i\in\mathcal{P}}(G_i-G_i^-)$. Thus the storage cost is reduced to $\bigO(n^2)$, which is scalable to the number of workers. Due to the space limitation, we do not provide further details in this work. 
	%Moreover, by invoking Sherman-Morrison-Woodbury matrix inverse formula \citep{horn2012matrix} we can also reduce the cost of computing $(G^{k})^{-1}$ from $\bigO(n^3)$ to $\bigO(n^2 \tau p)$ if $\tau p < n$. 
\end{rem}

\section{Numerical results}\label{sec:experiments}
In this section, we test the AGQN of Alg. \ref{alg:AGQN} and DAGQN of Alg. \ref{alg:DAGQN}  on the following logistic regression problem: 
\bea\label{equ:logistic_obj}
f(x) = \sum_{j=1}^{m}\ln(1 + e^{-b_{j}}\langle c_{j}, x \rangle) + \frac{\gamma}{2}\Vert x \Vert^{2}, \gamma > 0
\ena
where $c_{j}\in\Rn$ and $b_{j}\in\lbrace -1, 1 \rbrace$ are training samples, and $\gamma$ is the regularization parameter.

\subsection{Verify the perfomance of the AGQN}\label{subsec:6_1}
We firstly verify the performance of the AGQN  on the \textit{phishing} dataset \citep{chang2011libsvm} where $n = 68$ and $m = 11055$,
and compare it with the Nesterov accelerated gradient descent (NAGD) \citep{nesterov2003introductory} and the BFGS in \eqref{equ:BFGS_update}  with line search to satisfy Wolfe condition \citep{nocedal2006numerical}.
%, which takes the following form\citep{nesterov2003introductory}:
%\bea\label{equ:nes_acc_gd}
%& x^{k+1} = y^{k} - \frac{1}{\omega}\nabla f(y^{k})\\
%& y^{k+1} = x^{k + 1} - \frac{\sqrt{\omega} - \sqrt{\mu}}{\sqrt{\omega}+ \sqrt{\mu}}(x^{k+1} - x^{k}).
%\ena
We  set $\mu=0.02m$, $\omega = m$, $L = 0.04m$, $M=2$ and $\tau=6$ for the \textsf{gBroyd} in AGQN and tune parameters for the NAGD and BFGS separately. 

\begin{figure}
	\centering
	\includegraphics[width=.45\linewidth]{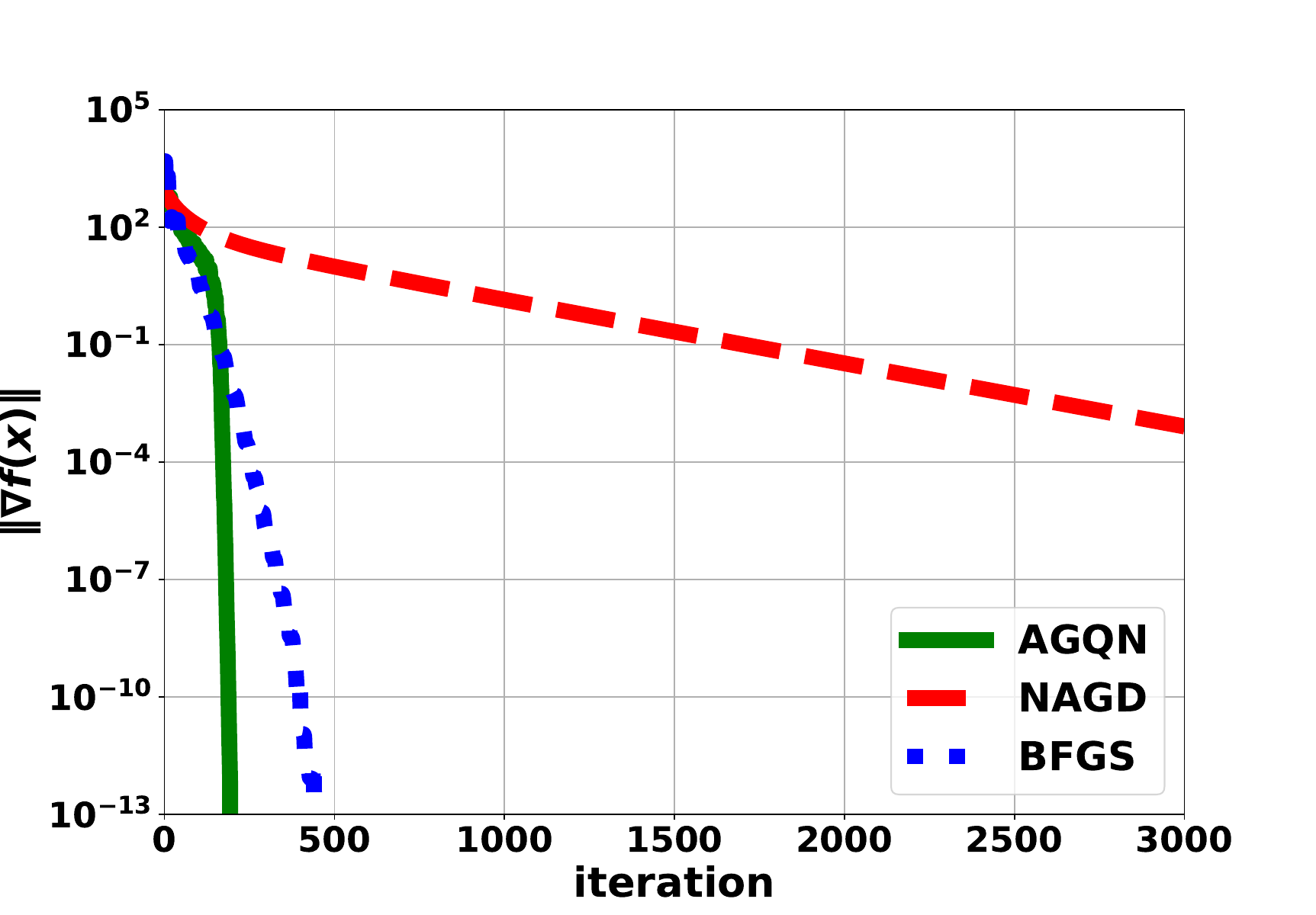}
	\caption{The convergence rates of NAGD, BFGS and the AGQN with respect to iteration number. }
	\label{fig:AGQN_fig1}
\end{figure}

 The initial point $x^{0}$ is randomly generated from the standard normal distribution.  Fig. \ref{fig:AGQN_fig1} depicts the convergence behavior of the algorithms in terms of iteration number.  As expected, the NAGD converges at a linear rate while the BFGS and AGQN achieve local superlinear convergence and the AGQN is essentially faster than BFGS.

\subsection{Verify the performance of the DAGQN}
In this subsection, we test the DAGQN  on the \textit{phishing} and the \textit{ijcnn1} datasets from LIBSVM repertory \citep{chang2011libsvm}, and compare with the distributed NAGD of first-order methods, and two second-order methods: DANE \citep{Shamir2014} and GIANT \citep{Wang2018}. Since both DANE and GIANT need to solve a sub-problem per iteration,  we adopt the conjugate gradient method \citep{nocedal2006numerical}. Besides, GIANT  requires a line search method on a  finite set of stepsizes for global convergence.

\begin{figure}
	\centering
	\subfigure[\textit{ijcnn1},  $\gamma = 1$]{\includegraphics[width=.45\linewidth]{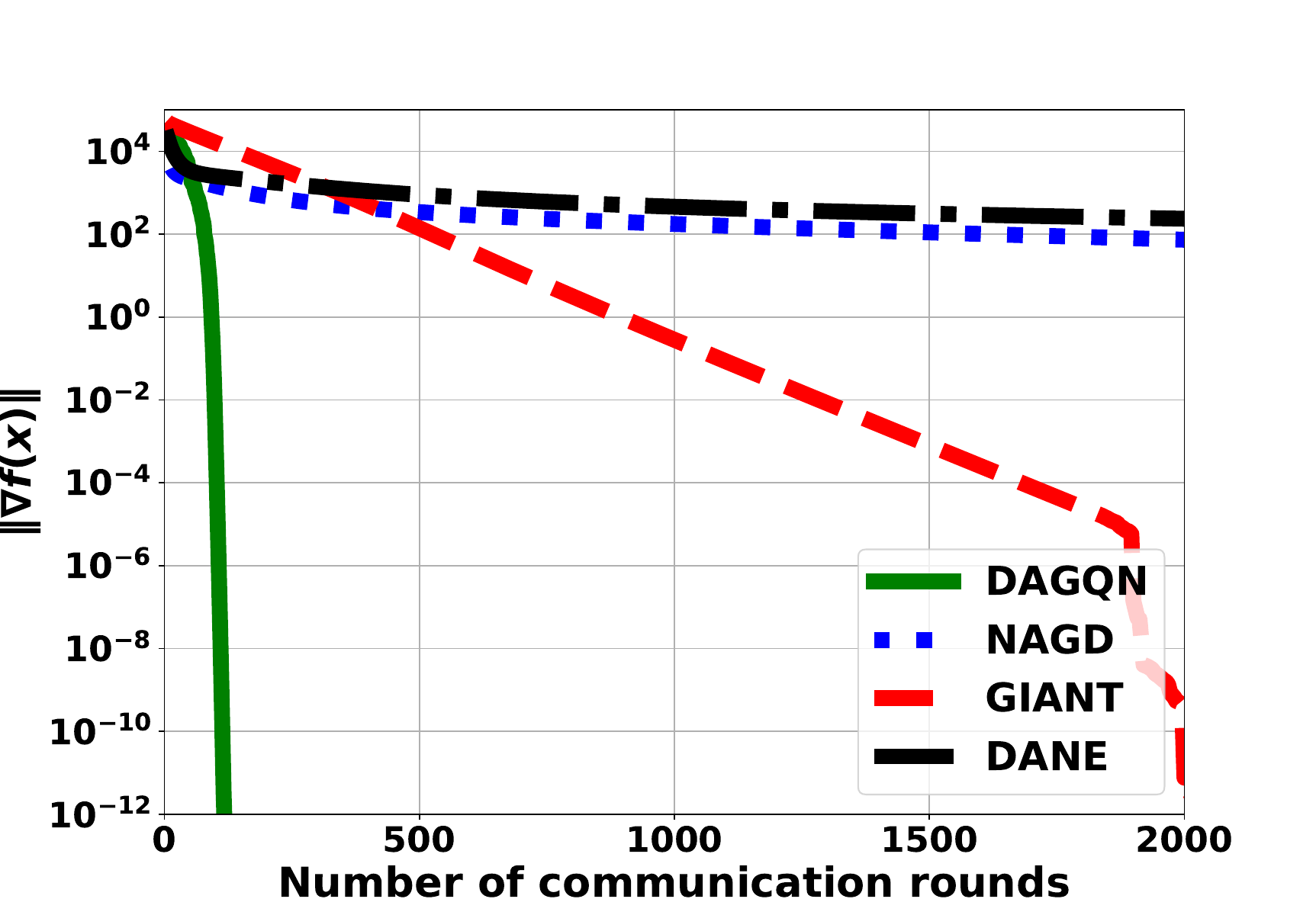}
		\label{fig:DAGQN_ijcnn1_comm_round}
	}
	\subfigure[\textit{ijcnn1}, $\gamma = 0.1$]{
		\includegraphics[width=.45\linewidth]{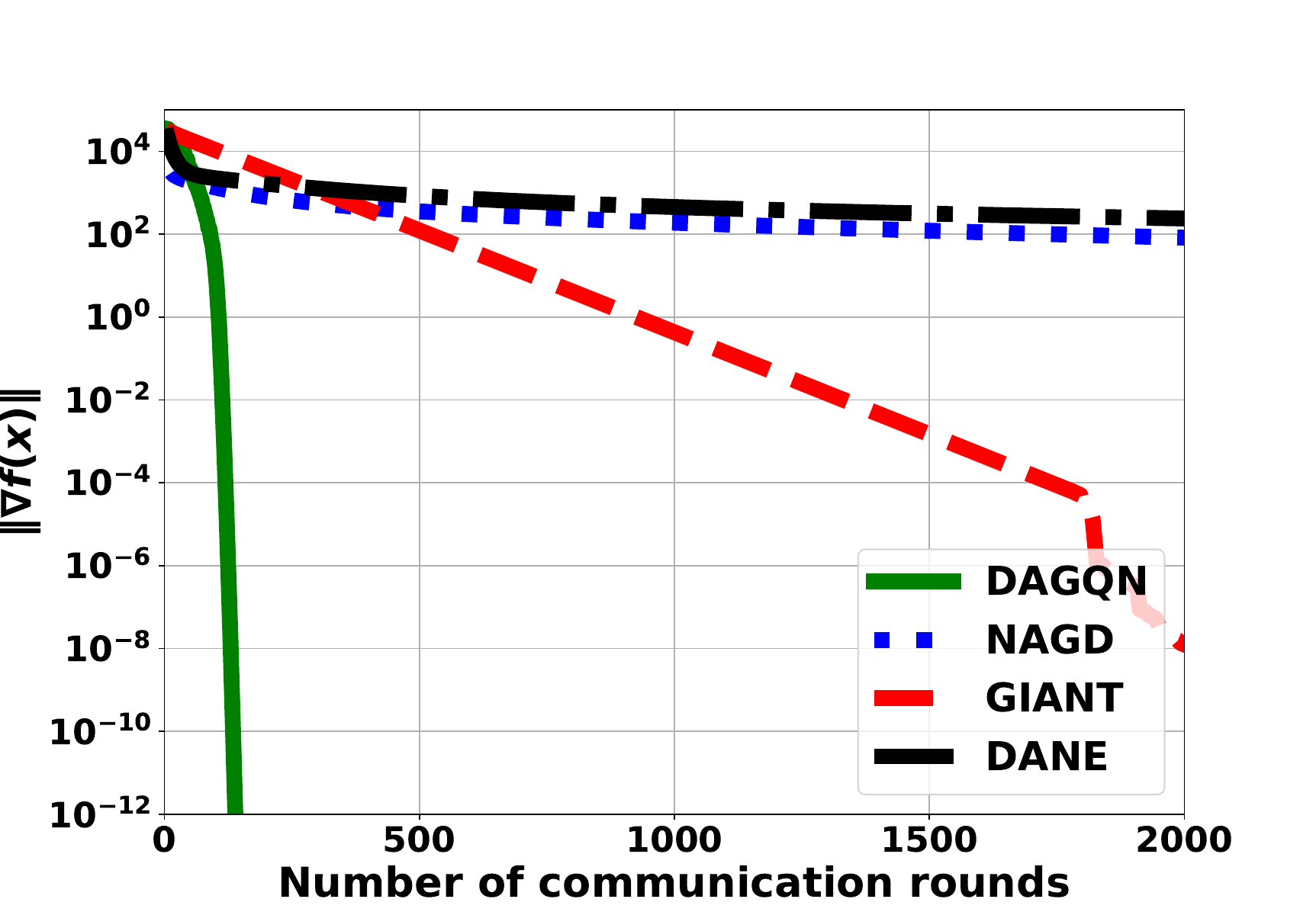}
		\label{fig:DAGQN_ijcnn1_gam01_comm_round}
	}
	\caption{$\|\nabla f\|$ versus the number of communication rounds between the master and each worker.}
			\label{fig:DAGQN_comm_round}
	\end{figure}

\begin{figure}
	\centering	
	\subfigure[\textit{ijcnn1}, $\gamma = 1$]{
		\includegraphics[width=.45\linewidth]{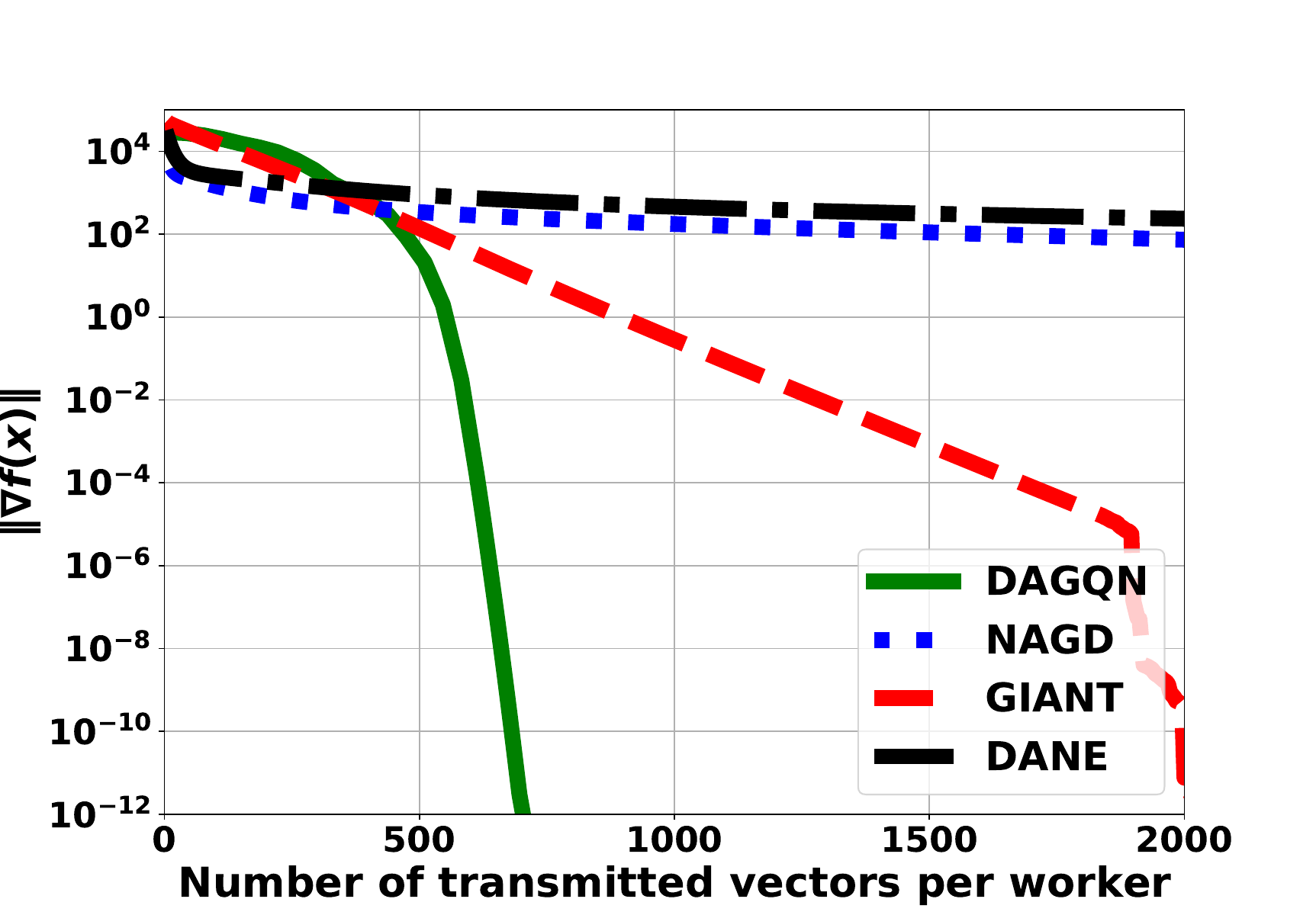}
		\label{fig:DAGQN_ijcnn1_comm_vecs}
	}
	\subfigure[ \textit{ijcnn1}, $\gamma = 0.1$]{
		\includegraphics[width=.45\linewidth]{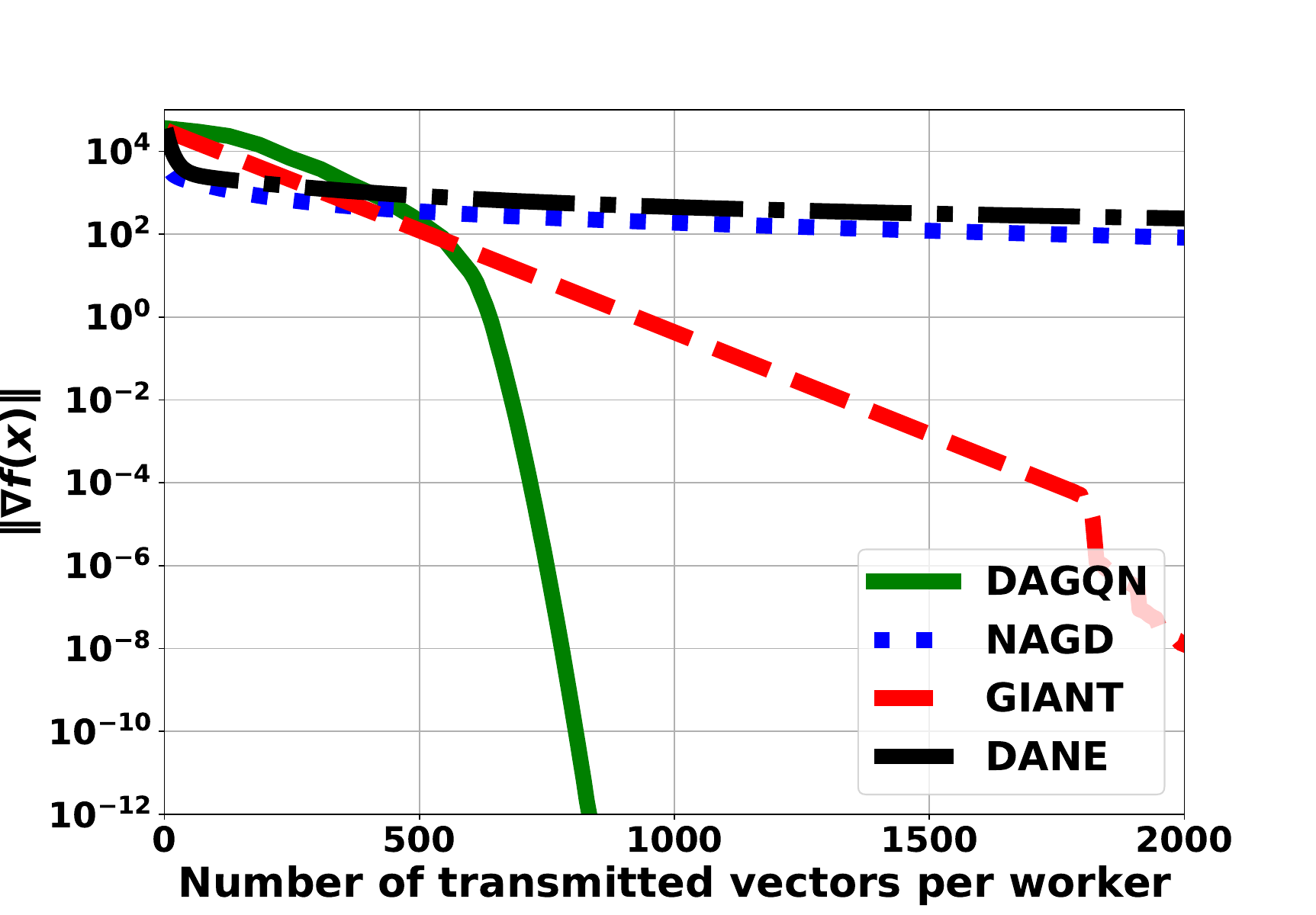}
		\label{fig:DAGQN_ijcnn1_gam01_comm_vecs}
	}
		\caption{$\|\nabla f\|$ versus the number of transmitted $n$-dimension vectors from each worker to the master.}
		\label{fig:DAGQN_comm_vecs}
	\end{figure}

 %For the \textit{phishing} dataset, we set $p=10$ workers and $\tau = 6$, and 
 Firstly, for the \textit{ijcnn1} dataset where $n = 22$ and $m = 141691$, we set $p=100$ workers and $\tau = 2$. All the training samples are divided equally among the workers. The parameters in DANE and GIANT are tuned via the methods in their original papers and are the same as that of Section \ref{subsec:6_1} for the DAGQN and NAGD, except that the total dataset size $m$ is replaced by $m / p$. 
 Now, we examine their convergence behaviors with respect to the communication cost between the master and a worker, which is measured by the number of communication rounds and the transmitted vectors. Then, we test two cases that $\gamma = 1$ and $\gamma=0.1$. Fig. \ref{fig:DAGQN_comm_round} confirms that the DAGQN significantly reduces the number of communication rounds. Note that both DANE and GIANT require multiple rounds of coordinated communications per iteration while  one round is sufficient for the DAGQN. A similar observation can also be found in Fig. \ref{fig:DAGQN_comm_vecs} for the case of transmitted $n$-dimension vectors. Besides, the smaller value of $\gamma = 0.1$ slightly degrades the performance of all these methods. But the DAGQN still works better than others.

\begin{figure}
	\centering	
	\includegraphics[width=.45\linewidth]{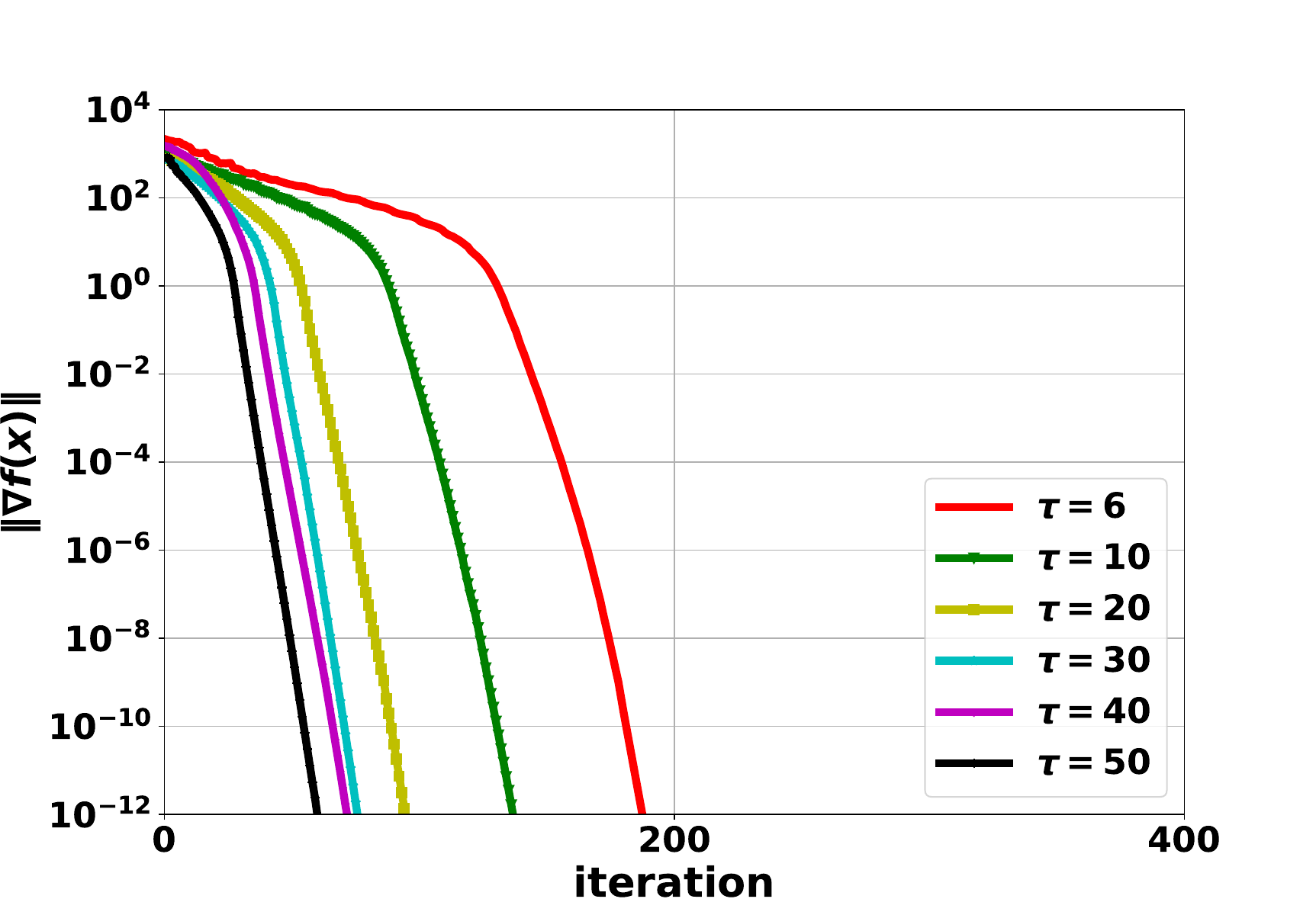}
	\caption{$\|\nabla f\|$ versus iteration number. }
	\label{fig:DAGQN_tau}
\end{figure}

Then, for the \textit{phishing} dataset, we set $p=10$ workers and $\tau = 6$. And we use Fig. \ref{fig:DAGQN_tau} to illustrate how the iteration number of DAGQN is affected by the parameter $\tau$, the repeated number of \textsf{gBroyd} updates per iteration. As expected, the total iteration number essentially decreases as $\tau$ increases.

\section{Conclusion}\label{sec:conclusion}
This paper has proposed two novel  adaptive quasi-Newton methods for the both the centralized and distributed settings. They all enjoy explicit non-asymptotic bounds for both the global convergence rate and local superlinear rate, though the communication complexity of the distributed one is comparable to that of the distributed first-order methods.  In the future, it is interesting to study their stochastic versions for the empirical risk problems.

\appendix
\section{Proof of Lemma \ref{lem:uniform_bound}}
\label{app:proof_uniform_bound}
\begin{pf}
	Letting the right-hand side of \eqref{equ:Broyd_inequ2} less than $\varepsilon$, we obtain a quadratic inequality in $r^+$, i.e., 
	\bea\label{equ:r_inequ}
	& \rho (\sigma_\tau + n) (Mr^+)^2 + 2\rho \sigma_{\tau}Mr^{+} + (\rho^{\tau + 1} - 1)\varepsilon \leq 0.
	\ena
	Since $\tau \geq 0$ and $\rho < 1$, the quadratic inequality indeed has a positive solution. Specifically, if $0 < Mr^+ \leq c_{\tau,\varepsilon}$ (see \eqref{equ:bar_r_tau}), then \eqref{equ:r_inequ} holds and $\sigma_{\nabla^{2}f(x^{+})}(G^+) \leq \varepsilon$. 
	
	%	Moreover, it is straightforward from \eqref{equ:r_inequ} that as $\tau$ gets larger, a larger $r^+$ is allowed to keep \eqref{equ:r_inequ} satisfied, which means that $c_\tau$ also gets larger.
	
	By \eqref{equ:x_update},  $r^+$ can be directly controlled by $\alpha$, i.e., 
	\bea
	r^+ &= \alpha \Vert G^{-1}\nabla f(x) \Vert_{\nabla^2 f(x)} \\
	& \leq \alpha \Vert G^{-1} \nabla f(x) \Vert_{G} = \alpha \Vert \nabla f(x) \Vert_{G}^{*}.
	\ena
	Plugging it into the above upper bound and solving the following inequality with respect to $\alpha\in(0,1]$
	\bea
	\alpha M \Vert \nabla f(x) \Vert_{G}^{*} \leq c_{\tau,\varepsilon},
	\ena
	we can easily obtain  \eqref{equ:bar_alpha_defi}. 
	
	Additionally, we give a concise lower bound for $c_{\tau, \varepsilon}$ to exploit its dependence on key constants such as $n$, $\tau$ and $\kappa$. Using \eqref{equ:bar_r_tau} and the ralation $\sqrt{1 + a} \geq 1 + \min\lbrace a/3, \sqrt{a/3} \rbrace$ for any $a \geq 0$, we have that
	\bea\label{equ:inequ_c}
	& c_{\tau, \varepsilon} \geq \\
	& \dfrac{\sigma_\tau}{\sigma_\tau + n} \min \left\lbrace \dfrac{\sigma_\tau + n}{3\rho\sigma_{\tau}^2} (1 - \rho^{\tau + 1}) \varepsilon, \sqrt{\dfrac{\sigma_\tau + n}{3\rho\sigma_{\tau}^2} (1 - \rho^{\tau + 1}) \varepsilon} \right\rbrace.
	\ena
	Applying the relation $(1 + a)^r\leq 1 / (1 - ra), a\in[-1, 1/r), r \geq 0$ to $\rho^{\tau + 1}$ (recall that $\rho = 1- 1(n\kappa)$), we have that
	\bea\label{equ:rho_tau_bound}
	 \rho^{\tau + 1} \leq \dfrac{n\kappa}{n\kappa + \tau + 1}.
	 \ena
	Substituting it into \eqref{equ:inequ_c}, we obtain the following lower bound expressed in $\Omega(\cdot)$
	\bea\label{equ:c_lower_bound}
	c_{\tau, \varepsilon} \geq \dfrac{ (\tau + 1) \varepsilon n \kappa}{ 3\rho ( 4(n\kappa)^2 \varepsilon + 4n(n\kappa)^2 ) } \geq \Omega \left( \dfrac{(\tau + 1) \varepsilon}{ (n + \varepsilon) (\kappa n + \tau) } \right).
	\ena
\end{pf}

\section{Proof of Lemma \ref{lem:adp_step}}
\label{app:proof_adp_step}
By Assumption \ref{assum:f} and Eq. \eqref{equ:G_error_inequ}, we have that  
\bea\label{equ:nabla_G_bounds}
& \mu I \preceq \nabla^2 f(x) \preceq \omega I\\
& \nabla^2 f(x) \preceq G \preceq (1 + \varepsilon)\nabla^2 f(x).
\ena
Then we bound $\Vert I - \nabla^2f(x)G^{-1}\alpha\Vert $ as follows
\bea 
& \Vert I - \nabla^2f(x)G^{-1}\alpha \Vert  \\
\
& = \Vert (1 - \alpha)I + \alpha \nabla^2f(x)(\nabla^2f(x)^{-1} - G^{-1})\Vert \\
\
& \leq (1 - \alpha) + \alpha \| \nabla^2f(x)  \| \| \nabla^2f(x)^{-1} - G^{-1} \|\\
\
& \leq (1 - \alpha) + \alpha \| \nabla^2f(x)  \| \| \nabla^2f(x)^{-1}\| \varepsilon / (1 + \varepsilon)\\
\
& \leq (1 - \alpha) + \alpha \varepsilon \kappa / (1 + \varepsilon) \leq 1- \alpha / 2
\ena
where the second inequality is due to relation \eqref{equ:nabla_G_bounds} and the last inequality comes from \eqref{equ:varepsilon_0}.

Now let $d = -G^{-1}\nabla f(x)$ and expand $\nabla f(x + \alpha d)$ as  
\bea\label{equ:deri1}
\nabla f(x + \alpha d) & = \nabla f(x) + (\nabla^{2}f(x)\cdot d)\alpha \\
& \quad + \left\lbrace \int_{0}^{\alpha} \left[  \nabla^{2}f(x + t d) - \nabla^2f(x) \right]dt\right\rbrace d\\
\
& = (I - \nabla^{2}f(x)G^{-1} \alpha)\nabla f(x)\\
& \quad + \left\lbrace \int_{0}^{\alpha} \left[  \nabla^{2}f(x + t d) - \nabla^2f(x) \right]dt\right \rbrace d.\\
\ena
Taking Euclidean norms on both sides leads to that
\bea
& \Vert \nabla f(x + \alpha d) \Vert \leq \Vert I - \nabla^{2}f(x)G^{-1} \alpha\Vert \Vert\nabla f(x) \Vert \\
& \quad\quad + \Vert d \Vert \left \lbrace \int_{0}^{\alpha} \Vert  \nabla^{2}f(x + t d) - \nabla^2f(x) \Vert dt \right\rbrace\\
&  \leq (1 - \alpha / 2 ) \Vert \nabla f(x) \Vert + L\alpha^{2}\Vert d \Vert^{2} / 2\\
&  \leq (1 - \alpha /2 ) \Vert \nabla f(x) \Vert + L \alpha^{2}\Vert \nabla f(x) \Vert^{2} / ( 2\mu^{2} ).
\ena
where Assumption \ref{assum:f}(b) is used in the first inequality and the last inequality is due to Assumption \ref{assum:f}(a) and \eqref{equ:G_condition}.

Multiplying both sides with $L / (2\mu^{2})$ yields that
\bea\label{equ:beta_inequ}
\beta(x^+) \leq (1 - \alpha / 2 )\beta(x) + \beta(x)^{2}\alpha^{2}.
\ena
Then, an ``optimal" stepsize can be derived by minimizing the quadratic function with respect to $\alpha$ on the right-hand side, the details of which are straightforward and thus omitted. 

\section{Proof of Theorem \ref{thm:AGQN_local}}\label{app:AGQN_local}
Firstly, we restate Lemma 4.3 in \citet{rodomanov2021greedy}, which is a handy tool to depict local evolution of $\lambda(x^k)$ when Hessian estimate error is bounded.
\begin{lem}\label{lem:lambda_local_evol}
	Let $x\in\R^n$, and let $G\in\R^{n\times n}$ be a positive definite matrix, satisfying
	\bea\label{equ:eta_condition}
	\nabla^2 f(x) \preceq G \preceq \eta \nabla^2 f(x)
	\ena
	for some $\eta \geq 1$. Let
	\bea
	x^+ \equdef x - G^{-1} \nabla f(x), 
	\ena
	and suppose that $M \lambda(x) \leq 2$. Then $r^+ \equdef \Vert x^+ - x \Vert_{\nabla^2 f(x)} \leq \lambda(x)$, and
	\bea
	\lambda(x^+) \leq \left(1 + \dfrac{ M\lambda(x)}{2} \right)\dfrac{\eta - 1 + M\lambda(x) / 2}{\eta} \lambda(x).
	\ena
\end{lem}

Condition \eqref{equ:AGQN_ll_cond1} ensures that the stepsize of the AGQN becomes 1. Jointly with Lemma \ref{lem:uniform_bound}, the initial condition $\sigma_{\nabla^2 f(x^0)} (G^0) \leq \varepsilon$ ensures that \eqref{equ:eta_condition} in Lemma \ref{lem:lambda_local_evol} holds with $\eta = 1 + \varepsilon_0$ for all $k$. Then, by applying Lemma \ref{lem:lambda_local_evol} for any $k$, we have that
\bea
\lambda(x^{k+1})  
& \leq \left(1 + \dfrac{M\lambda(x^k)}{2} \right) \dfrac{\varepsilon_0 + M\lambda(x^k) / 2}{1 + \varepsilon_0} \lambda(x^k)\\
\
& \leq \left( 1 - \dfrac{1}{1 + \varepsilon_0} + \dfrac{M\lambda(x^k)}{2} + \dfrac{ (M\lambda(x^k))^2 } { 4 (1 + \varepsilon_0) } \right) \lambda(x^k)\\
\
& \leq \left( \dfrac{1}{2} + \dfrac{M\lambda(x^k)}{2} + \dfrac{ ( M\lambda(x^k) )^2 } {4} \right) \lambda(x^k) \\
\ena
where the last inequality comes from the fact that $\varepsilon_0 \in (0, 1]$ (see \eqref{equ:varepsilon_0} in Lemma \ref{lem:adp_step}). Then, by condition \ref{equ:AGQN_ll_cond2} and simple induction, we have that $\lambda(x^k)$ is monotonically decreasing, 
\bea
M \lambda(x^{k+1}) \leq M \lambda(x^k) \leq \cdots \leq M\lambda(x^1) \leq M\lambda(x^0) \leq 1/4.
\ena
Thus, we have that
\bea
\lambda(x^{k+1})  
& \leq \left(\dfrac{1}{2} + \dfrac{M\lambda(x^k)}{2} + \dfrac{ ( M\lambda(x^k) )^2 } {4} \right) \lambda(x^k) \\
\
& \leq \left(\dfrac{1}{2} + M\lambda(x^k) \right) \lambda (x^k)  \leq \dfrac{3}{4} \lambda(x^k).
\ena

Now we establish the superlinear rate. 
For the simplicity of notation, let $\lambda_k \equdef \lambda ( x^k )$ and $\sigma_k \equdef \sigma_{\nabla^2 f(x^k)} (G^k) $. Define the following potential function:
\bea
\Phi_k \equdef \sigma_k + 4nM \lambda_k.
\ena
By invoking Lemma \ref{lem:lambda_local_evol} again and we have that
\bea\label{equ:poten_lambda}
\lambda_{k+1} 
& \leq \left(1 + \dfrac{ M \lambda_k }{2} \right) \dfrac{ \sigma_k + M \lambda_k / 2}{ 1 + \sigma_k} \lambda_k \\
\
& \leq \left(1 + \dfrac{ M \lambda_k }{2} \right) \left( \sigma_k + 4nM \lambda_k \right) \lambda_k \\
\
& \leq e^{ M \lambda_k / 2 } \Phi_k \lambda_k \leq e^{ 2 M \lambda_k } \Phi_k \lambda_k. 
\ena
Further, by Lemma 4.8 in \citet{rodomanov2021greedy} and the fact that our AGQN conducts $\tau$ additinonal \textsf{gBroyd} updates per iteration, we have that
\bea\label{equ:poten_sigma}
\sigma_{k+1} & \leq \rho \left(1 + M r^{k+1} \right)^2 \left( \rho^\tau \sigma_k + \dfrac{2nMr^{k+1} }{ 1 + Mr^{k+1} } \right) \\
\
& \leq \rho \left(1 + M \lambda_k \right)^2 \left( \rho^\tau \sigma_k + 2nM \lambda_k \right) \\
\
& \leq \rho e^{2M\lambda_k} \left( \rho^\tau \sigma_k + 2nM \lambda_k \right) \\
\
& \leq e^{2M\lambda_k} \rho^{\tau + 1} \sigma_k + \dfrac{3}{4} \rho e^{2M\lambda_k} \cdot 4nM \lambda_k\\
\
& \leq \max \lbrace \rho^{\tau + 1}, 3\rho / 4 \rbrace \cdot e^{ 2M\lambda_k } (\sigma_k + 4nM \lambda_k ) 
\\
\
& \leq \max \lbrace \rho^{\tau + 1}, 3 \rho / 4 \rbrace \cdot e^{ 2M\lambda_k } \Phi_k.\\
\ena
Then we focus on the evolution of $\Phi_k$:
\bea
& \Phi_{k+1} = \sigma_{k+1} + 4nM\lambda_{k+1} \\
\
& \leq \max \lbrace \rho^{\tau + 1}, 3 \rho / 4 \rbrace \cdot e^{ 2M\lambda_k } \Phi_k + 4nM e^{ 2 M \lambda_k } \Phi_k \lambda_k \\
\
& \leq \max \lbrace \rho^{\tau + 1}, 3 \rho / 4 \rbrace \cdot e^{ 2M\lambda_k } \Phi_k  \left( 1 + \dfrac{4nM\lambda_k}{ \max \lbrace \rho^{\tau + 1}, 3\rho / 4 \rbrace } \right).
\ena
To avoid the trivial case, we assume that $n \geq 2$, then we have that
\bea\label{equ:numb}
\max\lbrace \rho^{\tau + 1}, 3\rho / 4 \rbrace \geq 3\rho /4 & = 3/4 \times (1 - 1 / (n \kappa) ) \\
\
& \geq 3/4 \times 1/2 = 3/8. 
\ena
Subsequently, 
\bea
\Phi_{k+1} & \leq \max \lbrace \rho^{\tau + 1}, 3\rho / 4 \rbrace \cdot e^{ 2M\lambda_k } \Phi_k \left(1 + 12nM\lambda_k \right)\\
\
& \leq \max \lbrace \rho^{\tau + 1}, 3\rho / 4 \rbrace \cdot \Phi_k e^{ 2M( 6n + 1 ) \lambda_k } \\
\
& \leq \left(\max \lbrace \rho^{\tau + 1}, 3\rho / 4 \rbrace \right)^{k+1} \cdot \Phi_0 e^{ 2M( 6n + 1 ) \sum_{i=0}^k \lambda_i }.
\ena
Using the local linear rate of $\lambda_k$ in the above, we have that
\bea
e^{ 2M( 6n + 1 ) \sum_{i=0}^k \lambda_i } & \leq e^{2 ( 6n + 1 ) M\lambda_0 \sum_{i=0}^k (3/4)^i} \\
\
& \leq e^{8( 6n + 1 ) M\lambda_0} \leq 2
\ena
where the last inequality uses condition \eqref{equ:AGQN_ls_cond2}. Thus we establish the linear convergence of the potential $\Phi_k$:
\bea
\Phi_{k+1} \leq \left(\max \lbrace \rho^{\tau + 1}, 3\rho / 4 \rbrace \right)^{k+1} \cdot 2 \Phi_0
\ena
where
\bea
\Phi_0 = \sigma_0 + 4nM\lambda_0 \leq \varepsilon_0 + n\ln2/(12n + 2) = \bar{\Phi}_0. 
\ena
Finally, we establish the superlinear convergence of $\lambda_k$:
\bea
& \lambda_{k+1} \leq e^{2M\lambda_k} \Phi_k \lambda_k \leq e^{2M (6n + 1) \lambda_k} \Phi_k \lambda_k \\
\
& \leq e^{2M (6n + 1) \lambda_k} \left(\max \lbrace \rho^{\tau + 1}, 3\rho / 4 \rbrace \right)^{k} \\
& \quad \cdot \Phi_0 e^{ 2M(6n+1 ) \sum_{i=0}^{k-1} \lambda_i} \lambda_k \\
\
& \leq \left(\max \lbrace \rho^{\tau + 1}, 3\rho / 4 \rbrace \right)^{k} \cdot \Phi_0 e^{ 2M(6n+1 ) \sum_{i=0}^{k} \lambda_i} \lambda_k \\
\
& \leq \left(\max \lbrace \rho^{\tau + 1}, 3\rho / 4 \rbrace \right)^{k} \cdot 2\Phi_0 \lambda_k \\
\
& \leq \left(\max \lbrace \rho^{\tau + 1}, 3\rho / 4 \rbrace \right)^{k} \cdot 2 \bar{\Phi}_0 \lambda_k. \\
\ena

\section{Proof of Corollary \ref{cor:AGQN}}\label{app:cor_AGQN}
Starting from any $x^0\in\R^n$, we first consider the complexity of global stage, i.e., the number of iterations required before reaching the superlinear rate in Theorem \ref{thm:AGQN_local}. By Lemma \ref{lem:adp_step}, we have that 
\bea\label{equ:beta_descent}
& \beta(x^{k+1}) \leq
\begin{cases} 
	\beta(x^k)- 1 / 16, \quad\quad \text{if}~\beta(x^k) > 1 / (4 \bar{\alpha}_\tau(x^k) ), \\
	\
	(1 - \bar{\alpha}_\tau(x^k) / 4 )\beta(x^k), \quad \text{otherwise}.
\end{cases}
\ena
Telescoping the second case,
\bea
& \beta ( x^{k+1} ) \leq (1 - \min\lbrace \alpha_\tau(x^k) ,1 \rbrace / 4 ) \beta(x^k) \\
\
& \leq \max\lbrace (3/4) \beta(x^k), (1 - \alpha_\tau(x^k) /4) \beta(x^k) \rbrace\\
\
& = \max\lbrace (3/4) \beta(x^k), (1 - c_{\tau, \varepsilon_0} / (4M\Vert \nabla f(x^k) \Vert_{G^k}^{*} )) \beta(x^k) \rbrace \\
\
& \leq \max\lbrace (3/4) \beta(x^k), (1 - \sqrt{\mu}c_{\tau, \varepsilon_0} / (4M\Vert \nabla f(x^k) \Vert )) \beta(x^k) \rbrace \\
\
& = \max\lbrace (3/4) \beta(x^k), \beta(x^k) - \sqrt{\mu}Lc_{\tau, \varepsilon_0} / ( 8M\mu^2 )   \rbrace
\ena
where the second inequality uses the fact that $G^k \succeq \nabla^2 f(x^k) \succeq \mu I$.
Combining it with the first case, we have that
\bea
& \beta(x^{k+1}) \leq \\
& \quad \max \left\lbrace \frac{3}{4} \beta(x^k),  \beta(x^k) - \frac{ \sqrt{\mu}Lc_{\tau, \varepsilon_0} } { 8M \mu^2 },  \beta(x^k) - \frac{1}{16} \right \rbrace. 
\ena
Equivalently,
\bea
& \Vert \nabla f(x^{k+1}) \Vert \\
&  \leq \max \left\lbrace \frac{3}{4} \Vert \nabla f(x^k) \Vert ,  \Vert \nabla f(x^k) \Vert  - \min \left \lbrace \frac{ \sqrt{\mu}c_{\tau, \varepsilon_0} } { 4M }, \frac{\mu^2}{8L} \right \rbrace \right \rbrace\\
\
& \leq \max \left\lbrace \frac{3}{4} \Vert \nabla f(x^k) \Vert ,  \Vert \nabla f(x^k) \Vert  - \frac{ \mu^2}{4L} \min \lbrace c_{\tau, \varepsilon_0}, 1/2\rbrace \right \rbrace
\ena
where we use the relation $M \leq L / \mu^{3/2} $.
Then the complexity to reach condition \eqref{equ:AGQN_ll_cond1} is bounded by
\bea
& \bigO \left( \max \left\lbrace  \log\left( \frac{M \Vert \nabla f(x^0) \Vert }{\sqrt{\mu} c_{\tau, \varepsilon_0}} \right), \log\left( \frac{2L \Vert \nabla f(x^0) \Vert }{ \mu^2 } \right) \right\rbrace  \right. \\
\
& \quad\quad  \left. + \frac{4L \Vert \nabla f(x^0) \Vert }{ \mu^2 \min\lbrace c_{\tau, \varepsilon_0}, 1/2 \rbrace} \right).
\ena
Using $M\leq L/\mu^{3/2}$, \eqref{equ:c_lower_bound} and \eqref{equ:varepsilon_0} and hiding the $\log$ term, it becomes that
\bea\label{equ:AGQN_comp_part1}
\widetilde{\bigO} \left( \frac {nL (\kappa n + \tau) \kappa }{ \mu^2 \tau}  \Vert \nabla f(x^0) \Vert \right).
\ena
For condition \ref{equ:AGQN_ls_cond2}, considering that $\lambda(x) \leq \Vert \nabla f(x) \Vert/\sqrt{\mu} $, its complexity can be similarly bounded by \eqref{equ:AGQN_comp_part1} and we omit it.

For the superlinear stage, using \eqref{equ:rho_tau_bound}, we have that
\bea
\max \lbrace \rho^{\tau+1}, 3\rho / 4 \rbrace & \leq \rho \max\lbrace n\kappa / (n\kappa + \tau), 3/4 \rbrace \\
\
& \leq 1 - \max\lbrace n\kappa / (n\kappa + \tau), 3/4 \rbrace / n\kappa \\
\
& = 1 - \max\lbrace 1 / (n\kappa + \tau), 3 / (4 n\kappa) \rbrace. \\
\ena
Then we can simply convert the local superlinear rate in \eqref{equ:AGQN_ls_rate} to the second part of complexity in Corollary \ref{cor:AGQN} by the same way as in \citet{Safaryan2022}. 

\section{Proof of Lemma \ref{lem:f_prop}}\label{app:f_prop}
The proof part (a) and (b) are trivial and we omit them. By Lemma \ref{lem:self_concordance}, for any $x,y,z,w\in\R^n$, we have that
\bea
& \nabla^2 f(x) - \nabla^2 f(y) = \sum_{i=1}^p ( \nabla^2 f_i(x) - \nabla^2 f_i (y) ) \\
\
& \leq M \sum_{i=1}^p \Vert x - y \Vert_{\nabla^2 f_i(z)} \nabla^2 f_i(w) \\
\
& \leq M \nabla^2 f(w) \sum_{i=1}^p \Vert x - y \Vert_{\nabla^2 f_i (z)} \\
\
& \leq M \nabla^2 f(w) p \left( (1 / p) \sum_{i=1}^p \Vert x - y \Vert^2_{\nabla^2 f_i (z)} \right)^{1/2} \\
\
& = \sqrt{p}M \nabla^2 f(w) \Vert x - y \Vert_{\nabla^2 f(z)} 
\ena
where the second inequality uses $\nabla^2 f_i(w) \succ 0$, the third inequality comes from the power mean inequality $(1/p)\sum_{i=1}^p a_i \leq ( (1/p) \sum_{i=1}^p a_i^2) ^ {1/2}$ for any $a_i \geq 0$. The proof is finished.

\section{Proof of Lemma \ref{lem:uniform_bound_Delta}}
Similar with the proof of Lemma \ref{lem:uniform_bound} (see Appendix \ref{app:proof_uniform_bound}), using Lemma \ref{lem:new_GQN_update_linear_rate} and \ref{lem:Delta_evol}, we actually want the following to hold:
\bea
\rho(1 + Mr^{+})^{2}\left( \rho^\tau \varepsilon + \dfrac{ 2n\sqrt{p} Mr^{+} }{1 + Mr^{+}} \right) \leq \varepsilon.
\ena
Rearrange it and we obtain the following inequality in $r^+$,
\bea\label{equ:r_inequ_Delta}
& \rho (\sigma_\tau + n\sqrt{p}) (Mr^+)^2 + 2\rho \sigma_{\tau}Mr^{+} + (\rho^{\tau + 1} - 1)\varepsilon \leq 0
\ena
where $\sigma_\tau = \rho^\tau \varepsilon + n\sqrt{p}$. It is straightforward to verify that if $Mr^+ \leq c_{\tau,\varepsilon}$ (see \eqref{equ:c_tau_Delta}), \eqref{equ:r_inequ_Delta} holds and thus $\Delta_{\nabla^{2}f(x^{+})}(G^+) \leq \varepsilon$. 

Then by $G = \sum_{i=1}^p G^i \succeq \sum_{i=1}^p \nabla^2 f_i(x) = \nabla^2 f(x)$, we have that
\bea
r^+ &= \alpha \Vert G^{-1} \nabla f(x) \Vert_{\nabla^2 f(x)} \\
& \leq \alpha \Vert G^{-1} \nabla f(x) \Vert_{G} = \alpha \Vert \nabla f(x) \Vert_{G}^{*}.
\ena
Plugging it into $Mr^+ \leq c_{\tau,\varepsilon}$, we have that
\bea
\alpha M \Vert \nabla f(x) \Vert_{G}^{*} \leq c_{\tau,\varepsilon}.
\ena
Solve it with respect to $\alpha$ and obtain \eqref{equ:bar_alpha_defi_Delta}. 

\section{Proof of Theorem \ref{thm:DAGQN_theo}}
\label{app:proof_DAGQN_theo}
Firstly we give a lemma to show that $\Delta_{\nabla^2 f(x)}(G)$ is an upper bound of $\sigma_{ \nabla^2 f(x) }(G)$.
\begin{lem}\label{lem:two_sigma_comparison}
	Suppose that $G = \sum_{i=1}^p G_i$ and $G_{i} \succeq \nabla^{2}f_{i}(x) \succ 0$ for each $i\in\mathcal{P}$. It holds that 
	\bea\label{equ:inequ_Delta}
	\sigma_{\nabla^2 f(x)}(G) \leq \Delta_{\nabla^2 f(x) }(G).
	\ena
\end{lem}
\begin{pf}
	By the definition of $\sigma_{A}(G)$ in \eqref{equ:G_measure_defi}, we have that
	\bea\label{equ:barsigma_bound}
	\sigma_{\nabla^2 f(x)}(G) &= \langle \nabla^{2}f(x)^{-1}, G - \nabla^{2}f(x)\rangle \\
	\
	&= \sum_{i=1}^{p} \langle \nabla^{2}f(x)^{-1}, G_{i} - \nabla^{2}f_{i}(x)\rangle \\
	\
	& \leq  \sum_{i=1}^{p} \langle \nabla^{2}f_{i}(x)^{-1}, G_{i} - \nabla^{2}f_{i}(x)\rangle\\
	\ 
	& = \sum_{i=1}^{p} \sigma_{\nabla^2 f_i(x)}(G_{i}) = \Delta_{\nabla^2 f(x)}(G)
	\ena
	where the inequality comes from the fact that $\nabla^2 f(x) =\sum_{i=1}^p \nabla^2 f_i(x) \succ \nabla^2 f_i (x)$ since $\nabla^2 f_i (x) \succ 0$ for all $i$.
\end{pf}

The following lemma quantifies how $\Delta_{\nabla^2 f(x)}(G)$ evolves as $x$ is updated to $x^+$. 
\begin{lem}\label{lem:Delta_evol}
	Let $x\in\R^n$, and let $G \equdef \sum_{i=1}^p G_i$ where $G_i\in\R^{n\times n}$ is a positive definite matrix satisfying
$\nabla^2 f_i (x) \preceq G_i$ for all $i \in \mathcal{P}$. Let $x^+ \in \R^n$, $r_i^+ \equdef \Vert x^+ - x\Vert_{\nabla^2 f_i(x)}$ and $G^+ \equdef \sum_{i=1}^p G_i^+ $ where
	\bea
G_i^{+} =\gbd(\tG^+, \nabla^{2} f_i(x^{+})), ~~~ \tG_i^+ = (1 + Mr_i^+)G_i.         
	\ena
	Then we have that
	\bea
	& \Delta_{\nabla^2 f(x^+) }(G^{+}) \\
	& \quad\quad \leq \rho (1 + Mr^{+})^{2}\left(\Delta_{\nabla^2 f(x) }(G) + \dfrac{ 2n\sqrt{p} Mr^{+} }{1 + Mr^{+}} \right)
	\ena
	where $r^+ \equdef \Vert x^+ - x\Vert_{\nabla^2 f(x)}$ and $\rho \equdef 1 - 1/(n\kappa)$.
\end{lem}
\begin{pf}
	Applying Lemma 4.8 to each $G_i$, we have that
	\bea
	& \sigma_{\nabla^2 f_i (x^+)}(G_i^+) \\
	& \quad\quad \leq \rho (1 + Mr_i^+)^2 \left( \sigma_{\nabla^2 f_i (x)} (G_i) + \dfrac{2nMr_i^+}{1 + Mr_i^+} \right).
	\ena
	Taking sum over $i$, we have that
	\bea
	& \Delta_{\nabla^2 f(x^+)} (G^+) \\
	& \leq \rho \sum_{i=1}^p (1 + Mr_i^+) \left( \sigma_{\nabla^2 f_i (x)} (G_i) (1 + Mr_i^+) + 2nMr_i^+ \right) \\
	\
	& \leq \rho (1 + Mr^+)\left( \sum_{i=1}^p \sigma_{\nabla^2 f_i (x)} (G_i) (1 + Mr_i^+) \right. \\
	& \quad \left. + 2nM\sum_{i=1}^p r_i^+ \right) \\
	\
	& \leq \rho (1 + Mr^+) \left( (1 + Mr^+) \Delta_{\nabla^2 f(x)} (G)  + 2n\sqrt{p}M r^+ \right)\\
	\
	& = \rho (1 + Mr^+)^2 \left( \Delta_{\nabla^2 f(x)} (G)  + \dfrac{2n\sqrt{p}M r^+}{1 + M r^+} \right),
	\ena
	which the second and third equalities come from the relation $r^+ = (\sum_{i=1}^p (r_i^+)^2 )^{1/2} \geq r_i^+$ and the power mean equality $(1/p)\sum_{i=1}^p a_i \leq ((1/p) \sum_{i=1}^p a_i^2 )^{1/2}$.
	
	\subsection{Proof of part (a)} With the initial conditions in Theorem \ref{thm:DAGQN_theo}, applying Lemma \ref{lem:uniform_bound_Delta} and Lemma \ref{lem:two_sigma_comparison}, we have that
	\bea\label{equ:G_inequ_DAGQN}
	\nabla^2 f(x^k) \preceq G^k & \preceq (1 + \sigma_{x}(G^k)) \nabla^2 f(x^k) \\
	\
	& \preceq (1 + \Delta_{x}(G^k)) \nabla^2 f(x^k) \\
	\
	& \preceq (1 + \epsilon_0) \nabla^2 f(x^k), ~~~ \forall k.
	\ena 
	Note that $f$ has the same condition number $\kappa$ as each $f_i$. The remaining proof is the same as that of Lemma \ref{lem:adp_step} (see the Appendix \ref{app:proof_adp_step}) except that $\alpha_{\tau}(x)$ now is given by Lemma \ref{lem:uniform_bound_Delta} and $\beta(x)$ for $f$ in \eqref{equ:intro_obj} is defined in \eqref{per_metric_DAGQN}.
	\end{pf}
	
	\subsection{Proof of part (b)} With the initial condition given in \eqref{equ:DAGQN_ll_cond1} and the global convergence guarantee, it is straightforward to verify that $\alpha(x^k) = 1$ for all $k\geq 0$. Given \eqref{equ:G_inequ_DAGQN} and the initial condition \eqref{equ:DAGQN_ll_cond2}, we can apply Lemma  \ref{lem:lambda_local_evol} for any $k$, and have that
	\bea
	\lambda(x^{k+1})  
	& \leq \left(1 + \dfrac{M_f \lambda(x^k)}{2} \right) \dfrac{\varepsilon_0 + M_f \lambda(x^k) / 2}{1 + \varepsilon_0} \lambda(x^k).
	\ena
	The remaining proof is the same as that of Theorem \ref{thm:AGQN_local} (see Appendix \ref{app:AGQN_local}) and we omit it.
	
	\subsection{Proof of part (c)}
	For the simplicity of notation, let $\lambda_k \equdef \lambda ( x^k )$ and $\Delta_k \equdef \Delta_{\nabla^2 f(x^k)} (G^k) $. Define the following potential function:
	\bea
	\Phi_k \equdef \sigma_k + 4n\sqrt{p} M \lambda_k.
	\ena
	Since all the conditions of part (b) are also satisfied in part(c), we can invoke Lemma \ref{lem:lambda_local_evol} again and have that
	\bea\label{equ:poten_lambda_DAGQN}
	\lambda_{k+1} 
	& \leq \left(1 + \dfrac{ M_f \lambda_k }{2} \right) \dfrac{ \Delta_k + M_f \lambda_k / 2}{ 1 + \Delta_k} \lambda_k \\
	\
	& \leq \left(1 + \dfrac{ M_f \lambda_k }{2} \right) \left( \Delta_k + 4n \sqrt{p} M \lambda_k \right) \lambda_k \\
	\
	& \leq e^{ M_f \lambda_k / 2 } \Phi_k \lambda_k
        \leq e^{ 2 \sqrt{p}M \lambda_k } \Phi_k \lambda_k 
	\ena
	where the second and last inequality use the relation $M_f \leq \sqrt{p}M$ in Lemma \ref{lem:f_prop}.
	
	Further, by Lemma \ref{lem:new_GQN_update_linear_rate} and \ref{lem:Delta_evol}, we have that
	\bea\label{equ:poten_sigma_DAGQN}
	\Delta_{k+1} & \leq \rho \left(1 + M r^{k+1} \right)^2 \left( \rho^\tau \Delta_k + \dfrac{2n\sqrt{p}Mr^{k+1} }{ 1 + Mr^{k+1} } \right) \\
	\
	& \leq \rho \left(1 + M \lambda_k \right)^2 \left( \rho^\tau \Delta_k + 2n\sqrt{p} M \lambda_k \right) \\
	\
	& \leq \rho e^{2M\lambda_k} \left( \rho^\tau \Delta_k + 2n\sqrt{p}M \lambda_k \right) \\
	\
	& \leq e^{2M\lambda_k} \rho^{\tau + 1} \Delta_k + \dfrac{3}{4} \rho e^{2M\lambda_k} \cdot 4n\sqrt{p} M \lambda_k\\
	\
	& \leq \max \lbrace \rho^{\tau + 1}, 3\rho / 4 \rbrace \cdot e^{ 2M\lambda_k } (\Delta_k + 4n\sqrt{p}M \lambda_k ) 
	\\
	\
	& \leq \max \lbrace \rho^{\tau + 1}, 3 \rho / 4 \rbrace \cdot e^{ 2M\lambda_k } \Phi_k.\\
	\ena
	Then we focus on the evolution of $\Phi_k$:
	\bea
	& \Phi_{k+1} = \Delta_{k+1} + 4n\sqrt{p}M\lambda_{k+1} \\
	\
	& \leq \max \lbrace \rho^{\tau + 1}, 3 \rho / 4 \rbrace \cdot e^{ 2M\lambda_k } \Phi_k + 4n\sqrt{p}M e^{ 2 \sqrt{p} M \lambda_k } \Phi_k \lambda_k \\
	\
	& \leq \max \lbrace \rho^{\tau + 1}, 3 \rho / 4 \rbrace \cdot \Phi_k  \left( 1 + \dfrac{4n\sqrt{p}M\lambda_k}{ \max \lbrace \rho^{\tau + 1}, 3\rho / 4 \rbrace } \right) \cdot \\
	& \quad e^{ 2\sqrt{p}M\lambda_k }.
	\ena
	By \eqref{equ:numb}, we have that
	\bea
	\max\lbrace \rho^{\tau + 1}, 3\rho / 4 \rbrace \geq = 3/8. 
	\ena
	Then, 
	\bea
	\Phi_{k+1} & \leq \max \lbrace \rho^{\tau + 1}, 3\rho / 4 \rbrace \cdot e^{ 2\sqrt{p}M\lambda_k } \Phi_k \left(1 + 12n\sqrt{p}M\lambda_k \right)\\
	\
	& \leq \max \lbrace \rho^{\tau + 1}, 3\rho / 4 \rbrace \cdot \Phi_k e^{ 2\sqrt{p}M( 6n + 1 ) \lambda_k } \\
	\
	& \leq \left(\max \lbrace \rho^{\tau + 1}, 3\rho / 4 \rbrace \right)^{k+1} \cdot \Phi_0 e^{ 2\sqrt{p}M( 6n + 1 ) \sum_{i=0}^k \lambda_i }.
	\ena
	Using the local linear rate of $\lambda_k$ in part (b), we have that
	\bea\label{equ:e_bound_DAGQN}
	e^{ 2\sqrt{p}M( 6n + 1 ) \sum_{i=0}^k \lambda_i } & \leq e^{2 ( 6n + 1 ) \sqrt{p} M\lambda_0 \sum_{i=0}^k (3/4)^i} \\
	\
	& \leq e^{8( 6n + 1 ) \sqrt{p} M\lambda_0} \leq 2
	\ena
	where the last inequality uses condition \eqref{equ:DAGQN_ls_cond2}. Thus we establish the linear convergence of the potential $\Phi_k$:
	\bea
	\Phi_{k+1} \leq \left(\max \lbrace \rho^{\tau + 1}, 3\rho / 4 \rbrace \right)^{k+1} \cdot 2 \Phi_0
	\ena
	where
	\bea
	\Phi_0 = \Delta_0 + 4n\sqrt{p}M\lambda_0 \leq \varepsilon_0 + n\ln2/(12n + 2) = \bar{\Phi}_0. 
	\ena
	Finally, we establish the superlinear convergence of $\lambda_k$:
	\bea
	& \lambda_{k+1} \leq e^{2\sqrt{p}M\lambda_k} \Phi_k \lambda_k 
	   \leq e^{2\sqrt{p} M (6n + 1) \lambda_k} \Phi_k \lambda_k \\
	\
	& \leq e^{2\sqrt{p} M (6n + 1) \lambda_k} \left(\max \lbrace \rho^{\tau + 1}, 3\rho / 4 \rbrace \right)^{k} \cdot \lambda_k \\ 
	& \quad \cdot \Phi_0 e^{ 2\sqrt{p}M(6n+1 ) \sum_{i=0}^{k-1} \lambda_i}  \\
	\
	& \leq \left(\max \lbrace \rho^{\tau + 1}, 3\rho / 4 \rbrace \right)^{k} \cdot \lambda_k  \Phi_0 e^{ 2\sqrt{p}M(6n+1 ) \sum_{i=0}^{k} \lambda_i} \\
	\
	& \leq \left(\max \lbrace \rho^{\tau + 1}, 3\rho / 4 \rbrace \right)^{k} \cdot 2\Phi_0 \lambda_k \\
	\
	& \leq \left(\max \lbrace \rho^{\tau + 1}, 3\rho / 4 \rbrace \right)^{k} \cdot 2 \bar{\Phi}_0 \lambda_k \\
	\ena
	where the fifth inequality follows from \eqref{equ:e_bound_DAGQN}. The proof is finished.

\bibliography{autosam}  

@article{rodomanov2022rates,
	author = {Rodomanov, Anton and Nesterov, Yurii},
	journal = {Mathematical Programming},
	number = {1},
	pages = {159--190},
	publisher = {Springer},
	title = {Rates of superlinear convergence for classical quasi-Newton methods},
	volume = {194},
	year = {2022}}

@inproceedings{du2022adaptive,
	author = {Du, Yubo and You, Keyou},
	booktitle = {61st IEEE Conference on Decision and Control, Cancun, Mexico, USA},
	title = {Greedy quasi-Newton methods with explicit superlinear convergence},
	year = {2022}}

@article{jin2022non,
	author = {Jin, Qiujiang and Mokhtari, Aryan},
	journal = {Mathematical Programming},
	pages = {1--49},
	publisher = {Springer},
	title = {Non-asymptotic superlinear convergence of standard quasi-{Newton} methods},
	year = {2022}}

@article{rodomanov2021new,
	author = {Rodomanov, Anton and Nesterov, Yurii},
	journal = {Journal of optimization theory and applications},
	number = {3},
	pages = {744--769},
	publisher = {Springer},
	title = {New results on superlinear convergence of classical quasi-Newton methods},
	volume = {188},
	year = {2021}}

@article{rodomanov2021greedy,
	author = {Rodomanov, Anton and Nesterov, Yurii},
	journal = {SIAM Journal on Optimization},
	number = {1},
	pages = {785--811},
	title = {Greedy quasi-Newton methods with explicit superlinear convergence},
	volume = {31},
	year = {2021}}

@book{Boyd2004,
	author = {Boyd, Stephen and Vandenberghe, Lieven},
	publisher = {Cambridge University Press},
	title = {Convex Optimization},
	year = {2004}}

@article{polyak2020new,
	author = {Polyak, Boris and Tremba, Andrey},
	journal = {Optimization Methods and Software},
	number = {6},
	pages = {1272--1303},
	title = {New versions of Newton method: step-size choice, convergence domain and under-determined equations},
	volume = {35},
	year = {2020}}

@article{eisen2017decentralized,
	author = {Eisen, Mark and Mokhtari, Aryan and Ribeiro, Alejandro},
	journal = {IEEE Transactions on Signal Processing},
	number = {10},
	pages = {2613--2628},
	title = {Decentralized quasi-Newton methods},
	volume = {65},
	year = {2017}}

@inproceedings{Crane2019,
	author = {Crane, Rixon and Roosta, Fred},
	booktitle = {Advances in Neural Information Processing Systems 32},
	pages = {9498--9508},
	title = {DINGO: Distributed Newton-Type Method for Gradient-Norm Optimization},
	year = {2019}}

@inproceedings{Soori2020,
	author = {Soori, Saeed and Mishchenko, Konstantin and Mokhtari, Aryan and Dehnavi, Maryam Mehri and Gurbuzbalaban, Mert},
	booktitle = {International Conference on Artificial Intelligence and Statistics},
	pages = {1965--1976},
	title = {DAve-QN: A distributed averaged quasi-Newton method with local superlinear convergence rate},
	year = {2020}}

@inproceedings{Wang2018,
	author = {Wang, Shusen and Roosta, Fred and Xu, Peng and Mahoney, Michael W},
	booktitle = {Advances in Neural Information Processing Systems 31},
	pages = {2332--2342},
	title = {GIANT: Globally Improved Approximate Newton Method for Distributed Optimization},
	year = {2018}}

@inproceedings{Shamir2014,
	author = {Shamir, Ohad and Srebro, Nati and Zhang, Tong},
	booktitle = {International Conference on Machine Learning},
	pages = {1000--1008},
	title = {Communication-efficient distributed optimization using an approximate newton-type method},
	year = {2014}}

@book{nesterov2003introductory,
	author = {Nesterov, Yurii},
	publisher = {Springer Science \& Business Media},
	title = {Introductory lectures on convex optimization: A basic course},
	volume = {87},
	year = {2003}}

@article{chang2011libsvm,
	author = {Chang, Chih-Chung and Lin, Chih-Jen},
	journal = {ACM Transactions on Intelligent Systems and Technology},
	number = {3},
	pages = {1--27},
	title = {LIBSVM: a library for support vector machines},
	volume = {2},
	year = {2011}}

@inproceedings{recht2011hogwild,
	author = {Recht, Benjamin and Re, Christopher and Wright, Stephen and Feng Niu},
	booktitle = {Advances in Neural Information Processing Systems 24},
	pages = {693--701},
	title = {Hogwild: A Lock-Free Approach to Parallelizing Stochastic Gradient Descent},
	year = {2011}}

@article{liu2015an,
	author = {Liu, Ji and Wright, Stephen J and Re, Christopher and Bittorf, Victor and Sridhar, Srikrishna},
	journal = {Journal of Machine Learning Research},
	number = {1},
	pages = {285--322},
	title = {An asynchronous parallel stochastic coordinate descent algorithm},
	volume = {16},
	year = {2015}}

@inproceedings{lian2015asynchronous,
	author = {Lian, Xiangru and Huang, Yijun and Li, Yuncheng and Liu, Ji},
	booktitle = {Advances in Neural Information Processing Systems 28},
	pages = {2737--2745},
	title = {Asynchronous Parallel Stochastic Gradient for Nonconvex Optimization},
	year = {2015}}

@inproceedings{leblond2017asaga,
	author = {Leblond, R{\'e}mi and Pedregosa, Fabian and Lacoste-Julien, Simon},
	booktitle = {Artificial Intelligence and Statistics},
	pages = {46--54},
	title = {ASAGA: asynchronous parallel SAGA},
	year = {2017}}

@inproceedings{lee2018distributed,
	author = {Lee, Ching-pei and Lim, Cong Han and Wright, Stephen J},
	booktitle = {ACM SIGKDD International Conference on Knowledge Discovery and Data Mining},
	pages = {1646--1655},
	title = {A distributed quasi-newton algorithm for empirical risk minimization with nonsmooth regularization},
	year = {2018}}

@inproceedings{Islamov2021,
	author = {Islamov, Rustem and Qian, Xun and Richtarik, Peter},
	booktitle = {International Conference on Machine Learning},
	pages = {4617--4628},
	title = {Distributed Second Order Methods with Fast Rates and Compressed Communication},
	year = {2021}}

@article{zhang2022distributed,
	author = {Zhang, Jiaqi and You, Keyou and Ba{\c{s}}ar, Tamer},
	journal = {Automatica},
	pages = {110156},
	title = {Distributed adaptive Newton methods with global superlinear convergence},
	volume = {138},
	year = {2022}}

@article{Mokhtari2016,
	author = {Mokhtari, Aryan and Ling, Qing and Ribeiro, Alejandro},
	journal = {IEEE Transactions on Signal Processing},
	number = {1},
	pages = {146--161},
	title = {Network Newton distributed optimization methods},
	volume = {65},
	year = {2016}}

@book{nocedal2006numerical,
	author = {Nocedal, Jorge and Wright, Stephen},
	publisher = {Springer Science \& Business Media},
	title = {Numerical optimization},
	year = {2006}}

@article{pearlmutter1994fast,
	author = {Pearlmutter, Barak A},
	journal = {Neural computation},
	number = {1},
	pages = {147--160},
	title = {Fast exact multiplication by the Hessian},
	volume = {6},
	year = {1994}}

@article{Du2021,
	author = {Du, Yubo and You, Keyou},
	journal = {IEEE Transactions on Control of Network Systems},
	number = {3},
	pages = {1212-1224},
	title = {Asynchronous Stochastic Gradient Descent Over Decentralized Datasets},
	volume = {8},
	year = {2021}}

@book{horn2012matrix,
	author = {Horn, Roger A and Johnson, Charles R},
	publisher = {Cambridge University Press},
	title = {Matrix Analysis},
	year = {2012}}

@inproceedings{crane2020dino,
	author = {Crane, Rixon and Roosta, Fred},
	booktitle = {International Conference on Machine Learning},
	pages = {2174--2184},
	title = {{DINO}: distributed Newton-type optimization method},
	year = {2020}}

@article{mishchenko2021regularized,
	author = {Mishchenko, Konstantin},
	journal = {arXiv preprint arXiv:2112.02089},
	title = {Regularized Newton Method with Global $ O (1/k^2) $ Convergence},
	year = {2021}}

@inproceedings{Malitsky2020,
	author = {Malitsky, Yura and Mishchenko, Konstantin},
	booktitle = {International Conference on Machine Learning},
	pages = {6702--6712},
	title = {Adaptive Gradient Descent without Descent},
	year = {2020}}

@inproceedings{Safaryan2022,
	author = {Mher Safaryan and Rustem Islamov and Xun Qian and Peter Richt{\'{a}}rik},
	booktitle = {International Conference on Machine Learning},
	pages = {18959--19010},
	title = {FedNL: Making Newton-Type Methods Applicable to Federated Learning},
	year = {2022}}

@inproceedings{Li2021,
	author = {Li, Yichuan and Gong, Yonghai and Freris, Nikolaos M. and Voulgaris, Petros and Stipanovi{\'c}, Du{\v s}an},
	booktitle = {60th IEEE Conference on Decision and Control },
	pages = {1689-1694},
	title = {{BFGS-ADMM} for Large-Scale Distributed Optimization},
	year = {2021}}

@inproceedings{Li2022,
	author = {Li, Yichuan and Voulgaris, Petros G. and Freris, Nikolaos M.},
	booktitle = {IEEE International Conference on Acoustics, Speech and Signal Processing},
	pages = {4268-4272},
	title = {A Communication Efficient Quasi-{Newton} Method for Large-Scale Distributed Multi-Agent Optimization},
	year = {2022}}

@article{Tutunov2019,
	author = {Tutunov, Rasul and Bou-Ammar, Haitham and Jadbabaie, Ali},
	journal = {IEEE Transactions on Automatic Control},
	number = {10},
	pages = {3983-3994},
	title = {Distributed Newton Method for Large-Scale Consensus Optimization},
	volume = {64},
	year = {2019}}

@article{Mansoori2020,
	author = {Mansoori, Fatemeh and Wei, Ermin},
	journal = {IEEE Transactions on Automatic Control},
	number = {7},
	pages = {2769-2784},
	title = {A Fast Distributed Asynchronous Newton-Based Optimization Algorithm},
	volume = {65},
	year = {2020}}

@InProceedings{Alimisis2021,
  author    = {Alimisis, Foivos and Davies, Peter and Alistarh, Dan},
  title     = {Communication-efficient distributed optimization with quantized preconditioners},
  booktitle = {International Conference on Machine Learning},
  year      = {2021},
  pages     = {196--206},
}

@Article{Fabbro2022,
  author  = {Fabbro, Nicol{\`o} Dal and Dey, Subhrakanti and Rossi, Michele and Schenato, Luca},
  title   = {A Newton-type algorithm for federated learning based on incremental Hessian eigenvector sharing},
  journal = {arXiv preprint arXiv:2202.05800},
  year    = {2022},
}

@Article{Agafonov2022,
  author  = {Agafonov, Artem and Kamzolov, Dmitry and Tappenden, Rachael and Gasnikov, Alexander and Tak{\'a}{\v{c}}, Martin},
  title   = {FLECS: A Federated Learning Second-Order Framework via Compression and Sketching},
  journal = {arXiv preprint arXiv:2206.02009},
  year    = {2022},
}
\bibliographystyle{agsm}       

\end{document}